\documentclass{informs3}
\OneAndAHalfSpacedXII 

\usepackage{natbib}
 \bibpunct[, ]{(}{)}{,}{a}{}{,}%
 \def\bibfont{\small}%
\usepackage{graphicx,epstopdf}
\usepackage{enumitem}
\usepackage[algo2e,vlined,ruled]{algorithm2e}
\usepackage{subfig}
\newcommand{\R}{{\mathbb R}}
\newcommand{\cN}{\mathcal{N}}
\newcommand{\cC}{\mathcal{C}}
\newcommand{\cL}{\mathcal{L}}

\newcommand{\cF}{\mathcal{F}}
\newcommand{\cS}{\mathcal{S}}
\newcommand{\Exp}{\mathbb{E}}

\newcommand{\ceil}[1]{\lceil #1 \rceil}
\TheoremsNumberedThrough     

\EquationsNumberedThrough    

\begin{document}


 \RUNAUTHOR{Zhang et al.} 

\RUNTITLE{An Efficient Stochastic Augmented Lagrangian-Type Algorithm}

\TITLE{Solving  Stochastic Optimization with Expectation Constraints Efficiently by a Stochastic Augmented Lagrangian-Type Algorithm}

\ARTICLEAUTHORS{%
\AUTHOR{Liwei Zhang}
\AFF{School of Mathematical Sciences,
 Dalian University
of Technology, 116023 Dalian, China, \EMAIL{lwzhang@dlut.edu.cn}}
\AUTHOR{Yule Zhang}
\AFF{School of Science,
Dalian Martime University,
116085 Dalian, China, \EMAIL{ylzhang@dlmu.edu.cn}}

\AUTHOR{Jia Wu}
\AFF{School of Mathematical Sciences,
 Dalian University
of Technology, 116023 Dalian, China, \EMAIL{wujia@dlut.edu.cn}}

\AUTHOR{Xiantao Xiao}
\AFF{School of Mathematical Sciences,
 Dalian University
of Technology, 116023 Dalian, China, \EMAIL{xtxiao@dlut.edu.cn}}
} 

\ABSTRACT{%
This paper considers the problem of minimizing a convex expectation function  with a set of  inequality  convex expectation constraints. We propose a  stochastic augmented Lagrangian-type algorithm, namely the stochastic linearized proximal method of multipliers, to solve this convex stochastic optimization problem. This algorithm can be roughly viewed as a hybrid of stochastic approximation and the traditional proximal method of multipliers. Under mild conditions, we show that this algorithm exhibits  $O(K^{-1/2})$ expected convergence rates for both objective reduction  and constraint violation if parameters in the algorithm are properly chosen,  where $K$ denotes the number of iterations.  Moreover, we show that, with high probability, the algorithm has  $O(\log(K)K^{-1/2})$ constraint violation bound and   $O(\log^{3/2}(K)K^{-1/2})$ objective bound. Numerical results demonstrate that the proposed algorithm is efficient.
}%



\maketitle
\section{Introduction}\label{sec:intro}
In this paper, we consider the following stochastic optimization problem
\begin{equation}\label{eq:1}
\begin{array}{rl}
 \min\limits_{x \in \cC} & f(x):=\mathbb{E}[F(x,\xi)]\\[4pt]
{\rm s.t.} & g_i(x):=\mathbb{E}[G_i(x,\xi)] \leq 0,\ i=1,\ldots,p.\\
\end{array}
\end{equation}
Here, $\cC \subset \R^n$ is a nonempty bounded closed convex set, $\xi$ is a random vector whose probability distribution   is supported on  $\Xi \subseteq \R^q$, $F: \cC \times \Xi \rightarrow \R$ and $G_i:\cC \times \Xi \rightarrow \R$, $i=1,\ldots, p$.  Let $\Phi$ be the feasible set of problem (\ref{eq:1}) as
\begin{equation}\label{eq:Phi}
\Phi:=\left\{x\in \cC: g_i(x) \leq 0,\ i=1,\ldots,p\right\}.
\end{equation}
We assume that
$$
\mathbb{E}[F(x,\xi)]= \int_{\Xi} F(x,\xi)dP(\xi),\, \mathbb{E}[G_i(x,\xi)]= \int_{\Xi} G_i(x,\xi)dP(\xi),\ i=1,\ldots, p
$$
are well defined and finite valued for every $x\in \cC$. Moreover, we assume that the  functions $F(\cdot,\xi)$ and $G_i(\cdot,\xi)$ are continuous and convex on $\cC$ for almost every $\xi$. Hence,  the expectation  functions $f(\cdot)$ and $g_i(\cdot,\xi)$ are continuous and convex on $\cC$.
Problems in the form of (\ref{eq:1}) are standard in stochastic programming  \citep{RS2003,Romisch2003} and also arise frequently in many practical applications such as machine learning \citep{SN2005,TFZ2016} and finance \citep{RU2000,DR2003}.

A computational difficulty of solving  (\ref{eq:1}) is that expectation is a multidimensional integral and it cannot be computed with a high accuracy for large dimension $q$. In order to handle this issue, a popular approach is to use stochastic approximation (SA) technique which is based on the following general assumptions: (i)
it is possible to generate  i.i.d. sample $\xi^1,\xi^2,\ldots,$ of realizations of random vector $\xi$;
(ii) there is an oracle, which, for any point $(x,\xi)\in \cC \times \Xi$ returns stochastic subgradients $v_0(x,\xi),\ v_1(x,\xi),\ \ldots,\  v_p(x,\xi)$ of $F(x,\xi)$, $G_1(x,\xi),\ \ldots,\ G_p(x,\xi)$ such that
$        v_i(x)=\mathbb{E}[v_i(x,\xi)],\ i=0,1,\ldots,p $
        are well defined and are subgradients of $f(\cdot)$, $g_1(\cdot)$, $\ldots$, $g_p(\cdot)$ at $x$, respectively, i.e., $v_0(x)\in \partial f(x)$, $v_i(x) \in \partial g_i(x)$, $i=1,\ldots,p$.

Since the pioneering paper \citep{RM21}, due to low demand for computer memory and cheap computation cost at every iteration, SA type algorithms become widely used in stochastic optimization and machine learning, see, e.g. \cite{Pflug1996,BCN2018}.  If $f(\cdot)$ is twice continuously differentiable and strongly convex, in the classical analysis it is shown that  the SA algorithm exhibits asymptotically optimal rate of convergence $\mathbb{E}[f(x^k)-f^*]=O(k^{-1})$, where $x^k$ is $k$th iterate and $f^*$ is the optimal value.  An important improvement  developed by  \cite{Polyak1990} and   \cite{PJ1992} suggests that, larger stepsizes of SA algorithm can be adopted by  consequently averaging  the obtained iterates.   Moreover,  \cite{Lan2009}
  show that, without assuming smoothness and strong convexity, a properly modified SA method achieves the convergence  rate  $O(k^{-1/2})$ and remarkably outperforms the sample average approximation (SAA) approach for a certain class of convex stochastic problems. After the seminal work \citep{Lan2009}, there are many significant results appeared, even for nonconvex stochastic optimization problems, see \cite{BCN2018,Lan2020} and references cited therein.
 Among all mentioned works, the feasible  set is an abstract closed convex set and none of these SA algorithms are applicable to expectation constrained problems. The main reason is that the computation of  projection $\Pi_{\Phi}$ is quite easy only when $\Phi$ is of a simple structure.
 However, when $\Phi$ is defined by  (\ref{eq:Phi}), the computation is prohibitive.

As a first attempt for solving expectation constrained stochastic optimization problems with stochastic approximation technique, \cite{LanZ2016} introduce a cooperative stochastic approximation (CSA) algorithm for solving (\ref{eq:1}) with single expectation constraint ($p=1$),  which is a stochastic counterpart of  Polyak's subgradient method \citep{Polyak1967}.  The authors show that CSA exhibits the optimal $O(1/\sqrt{K})$ rate of expected convergence, where $K$ is a fixed iteration number. In an online fashion, \cite{YMNeely2017} propose an algorithm (simply denoted by ``YNW'') that can be easily extended to solve  (\ref{eq:1}) with multiple expectation constraints. Under the Slater's condition and the assumption that $\cC$ is compact, they show that the algorithm
can  achieve $O(1/\sqrt{K})$ expected regret  and $O(\log(K)/\sqrt{K})$
high probability regret. \cite{Xiao2019} develops a penalized stochastic gradient (PSG) method and establishes its almost sure convergence and expected convergence rates. PSG can be roughly viewed as a hybrid of the classical   penalty method  for nonlinear programming and the stochastic quasi-gradient method \citep{WFL2017} for stochastic composition problem.   A stochastic level-set method \citep{LNSY2019},  which ensures a feasible solution path with high probability, is proposed  and analyzed.
 \cite{ABR2021} propose a conservative stochastic optimization algorithm (CSOA), which is in the similar primal-dual framework as PSG and YNW. In addition to CSOA, the authors also propose a projection-free algorithm named as FW-CSOA which can deal with the case that the projection $\Pi_{\cC}$ is difficult to calculate. \cite{YX2022} study an adaptive primal-dual stochastic gradient method (APriD)  for solving (\ref{eq:1}) and establish the convergence rate of $O(1/\sqrt{K})$ in terms of the objective error and the constraint violation.

All of the above mentioned methods for solving (\ref{eq:1}) can be cast into the family of stochastic first-order algorithms. Although the iteration in stochastic first-order algorithms is computationally cheap and these methods  perform well for certain problems, there are  plenty of practical experiences and evidences of their convergence difficulties and instability with respect to the choice of parameters.
Recently,
 the success of augmented Lagrangian methods for various kinds of functional constrained optimization problems is witnessed. \cite{PR2007} study an augmented Lagrangian method for multistage stochastic problems. For solving semidefinite programming (SDP) problems, \cite{ZST2010} consider an Newton-CG augmented Lagrangian method,  which is shown to be very efficient even for large-scale SDP problems. \cite{DMW2016} propose several methods based on augmented Lagrangian framework for optimization problems with stochastic-order constraints and analyze their convergence. \cite{BSZ2021} study an augmented Lagrangian decomposition method for nonconvex chance-constrained problems, in which a convex subproblem and a 0-1 knapsack subproblem are solved at each iteration.
 The aim of this paper is to develop an efficient stochastic approximation-based augmented Lagrangian-type  method for solving (\ref{eq:1}). To the best of our knowledge, this is still limited in the literature.


 \cite{ZZW2020} propose a stochastic proximal method of multipliers (PMMSopt) for solving problem (\ref{eq:1}) and show that PMMSopt exhibits $O(K^{-1/2})$ convergence rates for both objective reduction and constraint violation. PMMSopt is partially inspired by the classic proximal method of multipliers  \citep{Rockafellar76B}, which is modeled through an augmented Lagrangian with an extra proximal term.
However, the subproblem is difficult to solve, that makes  PMMSopt an unimplementable algorithm, and hence no numerical results are given.

In this paper, based on PMMSopt,  we propose a stochastic linearized proximal method of multipliers (SLPMM) for solving the stochastic convex optimization  problem (\ref{eq:1}), and analyze its expected convergence rate as well as probability guarantee for both objective reduction and constraint violation.  In specific, at the $k$th iteration in SLPMM, we consider the augmented Lagrangian function $\cL_\sigma^k(x,\lambda)$ of a linearized problem with respect to the stochastic subgradients $v_i(x^k,\xi^k)$, $i=0,1,\ldots,p$. Then, we obtain the next iterate $x^{k+1}$ by solving the problem $\min_{x\in\cC}\cL_\sigma^k(x,\lambda^k)+\frac{\alpha}{2}\|x-x^k\|^2$ and update the Lagrange multiplier. The subproblem is  the minimization of a strongly convex (approximately) quadratic function   and hence is relatively easy to solve.   Assuming that the set $\cC$ is compact, the subgradients are bounded and the Slater's condition holds, if the parameters in SLPMM are chosen as $\alpha=\sqrt{K}$ and $\sigma=1/\sqrt{K}$, we show that SLPMM attains $O(1/\sqrt{K})$ expected convergence rate with respect to both objective reduction and constraint violation. Under certain light-tail assumptions, we also establish the large-deviation properties of SLPMM. The numerical results on some practical applications such as Neyman-Pearson classification  demonstrate that SLPMM performs efficiently and has certain advantages over the existing stochastic first-order methods.

%

The remaining parts of this paper are organized as follows. In Section \ref{sec:properties}, we develop some important properties of SLPMM. In Section \ref{sec:rates}, in the expectation sense we establish the convergence rate of SLPMM  for problem (\ref{eq:1}).  The high probability guarantees for objective reduction and constraint violation of SLPMM are investigated in Section \ref{sec:prob}. In  Section \ref{sec:num}, we  report our numerical results. Finally, we draw a conclusion  in Section \ref{sec:conclusion}.

\section{Algorithmic framework, assumptions and auxiliary lemmas}\label{sec:properties}
In this section, we propose a stochastic  linearized proximal method of multipliers (SLPMM) for solving  problem (\ref{eq:1}) and establish some important auxiliary lemmas.

Let us define $[t]_+:=\max\{t,0\}$ for any $t\in\R$ and let $[y]_+=\Pi_{\R^p_+}[y]$ denote  the projection of $y$ onto $\R^p_+$ for any $y \in \R^p$. We also define $[t]_+^2:=(\max\{t,0\})^2$.

The detail of SLPMM is described in Algorithm \ref{alg:SLPMM}.
\begin{algorithm2e}[htp]
\caption{A stochastic  linearized proximal method of multipliers}
\label{alg:SLPMM}
\lnlset{alg:SA1}{1}{Initialization: ~Choose an initial point $x^0 \in \cC$ and select parameters $\sigma>0,\alpha>0$.  Set $\lambda^0=0\in\R^p$ and $k=0$.}
\\ \vspace{.5ex}
\lnlset{alg:SA2}{2}\For{$k=0,1,2,\ldots$}{
\lnlset{alg:SA3}{3}{Generate i.i.d. sample $\xi^k$ of $\xi$ and compute
\be\label{xna}
x^{k+1}= \argmin\limits_{x \in \cC} \,\left\{ \cL^k_{\sigma }(x,\lambda^k) +\frac{\alpha}{2}\|x-x^k\|^2\right\},
\ee
where
\be\label{augL}
\begin{array}{ll}
\cL^k_{\sigma}(x,\lambda)&:=F(x^k,\xi^k)+\langle v_0(x^k,\xi^k),x-x^k \rangle\\[8pt]
& + \frac{1}{2\sigma}\left[ \displaystyle\sum_{i=1}^p[\lambda_i+\sigma (G_i(x^k,\xi^k)+ \langle v_i(x^k,\xi^k),x-x^k \rangle)]_+^2-\|\lambda\|^2\right]
\end{array}
\ee
and $v_i(x^k,\xi^k)$, $i=0,1,\ldots,p$ are the corresponding stochastic subgradients.
} \\ \vspace{.5ex}
\lnlset{alg:SA4}{4}{Update the Lagrange multipliers by
\be\label{xna1}
\lambda_i^{k+1}=[\lambda_i^k+\sigma (G_i(x^k,\xi^k)+ \langle v_i(x^k,\xi^k),x^{k+1}-x^k \rangle)]_+,\quad i=1,\ldots,p.
\ee
} \\ \vspace{.5ex}
\lnlset{alg:SA5}{5}{Set $k\leftarrow k+1$.}
\\ }
\end{algorithm2e}
In specific, at each iteration, we first generate an i.i.d. sample $\xi^k$ and choose the stochastic subgradients $v_i(x^k,\xi^k)$, $i=0,1,\ldots,p$ of $F$ and $G_i$, respectively.
Then, in (\ref{xna}) we obtain $x^{k+1}$ by computing the proximal point of $\cL^k_{\sigma}(x,\lambda)$, which is the augmented Lagrangian function of the linearized problem
\[
\begin{array}{ll}
 \min\limits_{x \in \cC} & F(x^k,\xi^k)+\langle v_0(x^k,\xi^k),x-x^k \rangle\\[4pt]
{\rm s.t.} & G_i(x^k,\xi^k)+ \langle v_i(x^k,\xi^k),x-x^k \rangle \leq 0,\ i=1,\ldots,p.\\
\end{array}
\]
Finally, in (\ref{xna1}) we update the Lagrange multipliers.

 Denote
\[
G(x,\xi):=(G_1(x,\xi),\ldots, G_p(x,\xi))^T,\quad g(x):=(g_1(x),\ldots, g_p(x))^T.
\]
 Let
\[
V(x^k,\xi^k):=(v_1(x^k,\xi^k),\ldots, v_p(x^k,\xi^k))^T,
\]
then (\ref{xna1}) can be rewritten as
\be\label{eq:xna1}
\lambda^{k+1}=[\lambda^k+\sigma (G(x^k,\xi^k)+V(x^k,\xi^k)(x^{k+1}-x^k))]_+.
\ee

In the following, we shall study the convergence of the stochastic process $\{x^k,\lambda^k\}$ generated by SLPMM with respect to the  filtration  $\cF_k$ (sigma-algebra) which is generated by the random information $\{(\xi^0,\ldots,\xi^{k-1})\}$.
Before that, we introduce the following assumptions.

\begin{assumption}\label{assu:compact} Let  $R>0$  be a positive parameter such that
        $$
        \|x'-x''\|\leq R,\ \forall x',x'' \in \cC.
        $$
\end{assumption}
\begin{assumption}\label{assu:cons}  There exists a constant $\nu_g>0$  such that for each $\xi^k$,
        $$
        \|G(x,\xi^k)\| \leq \nu_g,\ \forall x \in \cC.
        $$
\end{assumption}
\begin{assumption}\label{assu:moment} There exist constants $\kappa_f>0$ and $\kappa_g>0$ such that for each $\xi^k$,
                $$
        \|v_0(x,\xi^k)\| \leq \kappa_f, \,\, \|v_i(x,\xi^k)\| \leq \kappa_g,\ i=1,\ldots, p,\ \forall x\in \cC.
        $$
\end{assumption}
\begin{assumption}\label{assu:slater}
The Slater's condition holds, i.e., there exist $\varepsilon_0>0$ and $\widehat x \in \cC$ such that
        $$
        g_i(\widehat x) \leq -\varepsilon_0, \,\, i =1,\ldots, p.
        $$
\end{assumption}
Assumption \ref{assu:compact} shows that $\cC$ is a compact convex set with diameter $R$. Assumption \ref{assu:cons} indicates that the constraint functions $G_i(\cdot,\xi^k)$ are bounded over $\cC$. This assumption is a bit restrictive, but it is required in the analysis of almost all existing stochastic methods for solving problem (\ref{eq:1}) \citep{LanZ2016,YMNeely2017,LNSY2019,Xiao2019}.  Assumption \ref{assu:moment} requires that the stochastic subgradients $v_i(\cdot,\xi^k)$ are bounded over $\cC$. Assumption \ref{assu:slater} is a standard Slater's condition for optimization problem with functional constraints.

The following auxiliary lemma will be used several times in the sequel.
\begin{lemma}\label{lem:opt-x}
 For any $z\in \cC$, we have
\begin{equation}\label{eq:opt-x-1}
\begin{array}{ll}
\displaystyle\langle v_0(x^k,\xi^k),x^{k+1}-x^k \rangle + \frac{1}{2\sigma}\|\lambda^{k+1}\|^2
+ \frac{\alpha}{2} \|x^{k+1}-x^k\|^2  \\[5pt]
\leq \displaystyle\langle v_0(x^k,\xi^k),z-x^k \rangle + \frac{1}{2\sigma}\left[ \sum_{i=1}^p[\lambda^k_i+\sigma (G_i(x^k,\xi^k)+\langle v_i(x^k,\xi^k), z-x^k \rangle)]_+^2\right]\\[15pt]
\quad\quad+ \displaystyle\frac{\alpha}{2}(\|z-x^k\|^2-\|z-x^{k+1}\|^2).
\end{array}
\end{equation}
In particular, if we take $z=x^k$, it yields
\begin{equation}\label{eq:opt-x-2}
\begin{array}{ll}
\displaystyle\langle v_0(x^k,\xi^k),x^{k+1}-x^k \rangle + \frac{1}{2\sigma}\|\lambda^{k+1}\|^2
+ \alpha \|x^{k+1}-x^k\|^2\\[10pt]
\leq  \displaystyle\frac{1}{2\sigma}\left[ \sum_{i=1}^p[\lambda^k_i+\sigma G_i(x^k,\xi^k)]_+^2\right].
\end{array}
\end{equation}
\end{lemma}
\proof{Proof.}
By using the optimality conditions, we have from (\ref{xna}) that $x^{k+1}$ satisfies
\begin{equation}\label{eq:aux-opt}
0\in \nabla_x\cL^k_{\sigma}(x^{k+1},\lambda^k)+\alpha (x^{k+1}-x^k)+\cN_{\cC}(x^{k+1}),
\end{equation}
where $\cN_{\cC}(x^{k+1})$ denotes the normal cone of $\cC$ at $x^{k+1}$ and
\[
\nabla_x\cL^k_{\sigma}(x^{k+1},\lambda^k)=v_0(x^k,\xi^k)+\sum_{i=1}^pv_i(x^k,\xi^k)\cdot[\lambda_i^k+\sigma (G_i(x^k,\xi^k)+ \langle v_i(x^k,\xi^k),x^{k+1}-x^k \rangle)]_+.
\]
Let us now consider the following auxiliary problem
\begin{equation}\label{eq:auxP}
\begin{array}{ll}
\min\limits_{x \in \cC} \,\langle v_0(x^k,\xi^k),x-x^k \rangle+ \frac{1}{2\sigma}\left[ \displaystyle\sum_{i=1}^p[\lambda_i^k+\sigma (G_i(x^k,\xi^k)+ \langle v_i(x^k,\xi^k),x-x^k \rangle)]_+^2\right]\\[15pt]
\quad\quad\quad +\frac{\alpha}{2}(\|x-x^k\|^2-\|x-x^{k+1}\|^2).
\end{array}
\ee
We can easily check that (\ref{eq:auxP}) is a convex optimization problem. Therefore, $\hat{x}$ is an optimal solution to (\ref{eq:auxP}) if and only if
\[
\begin{array}{ll}
0\in &v_0(x^k,\xi^k)+\sum_{i=1}^pv_i(x^k,\xi^k)\cdot[\lambda_i^k+\sigma (G_i(x^k,\xi^k)+ \langle v_i(x^k,\xi^k),\hat{x}-x^k \rangle)]_+\\[6pt]&\quad+\alpha (x^{k+1}-x^k)+\cN_{\cC}(\hat{x}).
\end{array}
\]
Hence, if follows from (\ref{eq:aux-opt}) that $x^{k+1}$ is an optimal solution to (\ref{eq:auxP}), which gives (\ref{eq:opt-x-1}) and (\ref{eq:opt-x-2}) obviously.
\Halmos\endproof

In what follows, we  estimate an upper bound  of $\|x^{k+1}-x^k\|$.
\begin{lemma}\label{lem:aux3}
Let  Assumptions \ref{assu:compact}-\ref{assu:moment} be satisfied. Then, if the parameters  satisfy $2\alpha-p\kappa_g^2\sigma>0$, we have
\[
\|x^{k+1}-x^k\|\leq  \frac{1}{\alpha
}(\kappa_f+ \sqrt{p}\kappa_g\|\lambda^k\|+
\sqrt{p}\nu_g\kappa_g\sigma).
\]
\end{lemma}
\proof{Proof.}  From  (\ref{eq:opt-x-2}) and Assumption \ref{assu:moment}, we have
\[
\alpha\|x^{k+1}-x^k\|^2\leq  \kappa_f\|x^{k+1}-x^k\|+
\frac{1}{2\sigma}\sum_{i=1}^p\left([a_i]_+^2-[b_i]_+^2\right),
\]
in which, for simplicity, we use
\[
a_i:=\lambda_i^k+\sigma G_i(x^k,\xi^k),\quad b_i:=\lambda_i^k+\sigma (G_i(x^k,\xi^k)+\langle v_i(x^k,\xi^k),x^{k+1}-x^k\rangle).
\]
Noticing that
\[
\begin{array}{ll}
[a_i]_+^2-[b_i]_+^2=([a_i]_++[b_i]_+)([a_i]_+-[b_i]_+)\\[8pt]
\leq (|a_i|+|b_i|)\cdot|a_i-b_i|\\[8pt]
\leq (2|a_i|+|b_i-a_i|)\cdot|a_i-b_i|\\[8pt]
=2|a_i|\cdot|a_i-b_i|+(a_i-b_i)^2\\[8pt]\
\leq 2|\lambda_i^k+\sigma G_i(x^k,\xi^k)|\cdot\sigma\kappa_g\|x^{k+1}-x^k\|+\sigma^2\kappa_g^2\|x^{k+1}-x^k\|^2,
\end{array}
\]
we obtain
\[
2\alpha\|x^{k+1}-x^k\|\leq 2\kappa_f+\sum_{i=1}^p(2\kappa_g|\lambda_i^k+\sigma G_i(x^k,\xi^k)|+\sigma\kappa_g^2\|x^{k+1}-x^k\|).
\]
If $2\alpha-p\kappa_g^2\sigma>0$, it yields
\[
\|x^{k+1}-x^k\|\leq \frac{2}{2\alpha-p\kappa_g^2\sigma}\left(\kappa_f+\sum_{i=1}^p(\kappa_g|\lambda_i^k+\sigma G_i(x^k,\xi^k)|\right).
\]
Therefore,  from the facts that $\sum_{i=1}^p|\lambda_i^k|\leq\sqrt{p}\|\lambda^k\|$ and \[\sum_{i=1}^p|G_i(x^k,\xi^k)|\leq\sqrt{p}\|G(x^k,\xi^k)\|\leq\sqrt{p}\nu_g,\]
the claim is obtained.
\Halmos\endproof

Under the Slater's condition, we derive the following conditional expected estimate of the multipliers.
\begin{lemma}\label{lem:aux4}
Let Assumption \ref{assu:slater} be satisfied. Then, for any $t_2 \leq t_1-1$ where $t_1$ and $t_2$ are positive integers,
\[
 \mathbb{E}\left[\langle \lambda^{t_1}, G(\widehat x, \xi^{t_1}) \rangle \,|\, \cF_{t_2}\right]
\leq -\varepsilon_0 \mathbb{E}\left[\|\lambda^{t_1}\| \,|\, \cF_{t_2}\right].
\]
\end{lemma}
\proof{Proof.}
For any $i \in \{1,\ldots,p\}$, noticing that $\lambda^{t_1}_i\in \cF_{t_1}$  and $\cF_{t_2}\subseteq \cF_{t_1}$ for $t_2 \leq t_1-1$, we have
\[
\begin{array}{ll}
\mathbb{E}\left[\lambda^{t_1}_iG_i(\widehat x, \xi^{t_1}) \,|\, \cF_{t_2}\right] &=\mathbb{E} \left[\mathbb{E}\left[\lambda^{t_1}_iG_i(\widehat x, \xi^{t_1}) \,|\, \cF_{t_1}\right]\,|\,\cF_{t_2}\right]\\[10pt]
&=\mathbb{E}\left[\lambda^{t_1}_ig_i(\widehat{x})\,|\,\cF_{t_2}\right]\\[10pt]
&\leq -\varepsilon_0 \mathbb{E}\left[\lambda^{t_1}_i \,|\, \cF_{t_2}\right].
\end{array}
\]
Summing the above inequality over $i \in \{1,\ldots,p\}$ yields
$$
 \mathbb{E}\left[\langle \lambda^{t_1}, G(\widehat x, \xi^{t_1}) \rangle \,|\, \cF_{t_2}\right]
\leq
-\varepsilon_0 \mathbb{E} \left[ \sum_{i=1}^p \lambda^{t_1}_i \,|\, \cF_{t_2}\right]
\leq -\varepsilon_0 \mathbb{E}\left[\|\lambda^{t_1}\| \,|\, \cF_{t_2}\right],
$$
by using  $\sum_{i=1}^p \lambda^{t_1}_i\geq \|\lambda^{t_1}\|$.
\Halmos\endproof

We next present some important relations of $\|\lambda^k\|$.
\begin{lemma}\label{lem:aux5}
Let Assumptions \ref{assu:compact}--\ref{assu:slater} be satisfied and $s > 0$ be an arbitrary integer.   Define  $\beta_0:=\nu_g+\sqrt{p}\kappa_gR$ and
 \begin{equation}\label{eq:theta9}
 \vartheta (\sigma,\alpha,s):= \frac{\varepsilon_0\sigma s}{2}+\sigma \beta_0(s-1)+ \frac{\alpha R^2}{\varepsilon_0s}+\frac{2\kappa_f R}{\varepsilon_0}+ \frac{\sigma \nu_g^2}{\varepsilon_0}.
 \end{equation}
Then, the following  holds:
\begin{equation}\label{eq:6}
|\|\lambda^{k+1}\|-\|\lambda^k\||\leq \sigma \beta_0
\end{equation}
and
\begin{equation}\label{eq:7}
\mathbb{E}\left [ \|\lambda^{k+s}\|-\|\lambda^k\| \,|\, \cF_k\right]
\leq \left
\{
\begin{array}{ll}
s \sigma \beta_0, & \mbox{if } \|\lambda^k\| < \vartheta (\sigma,\alpha, s),\\[6pt]
-s \displaystyle \frac{\sigma \varepsilon_0}{2}, & \mbox{if } \|\lambda^k\| \geq  \vartheta (\sigma,\alpha,s).
\end{array}
\right.
\end{equation}

\end{lemma}
\proof{Proof.}
It follows from Assumptions \ref{assu:compact}--\ref{assu:moment}, (\ref{eq:xna1}) and the nonexpansion property of the projection $\Pi_{\R^p_+}(\cdot)$ that
\[
\begin{array}{ll}
&|\|\lambda^{k+1}\|-\|\lambda^k\||\\[6pt]
&\leq\|\lambda^{k+1}-\lambda^k\|  =\|[\lambda^k+\sigma (G(x^k,\xi^k)+V(x^k,\xi^k)(x^{k+1}-x^k))]_+-[\lambda^k]_+\|\\[6pt]
& \leq
\sigma \|G(x^k,\xi^k)+V(x^k,\xi^k)(x^{k+1}-x^k)\|\\[6pt]
& \leq \sigma [\nu_g+\sqrt{p}\kappa_g R],
\end{array}
\]
which implies (\ref{eq:6}).  This also gives that $\|\lambda^{k+s}\|-\|\lambda^k\|\leq s \sigma \beta_0$. Hence, we only need to establish the second part  in (\ref{eq:7}) under the case $\|\lambda^k\| \geq  \vartheta (\sigma,\alpha,s)$.

For a given positive integer $s$, suppose that $\|\lambda^k\| \geq \vartheta (\sigma,\alpha,s)$. For any $l \in \{k,k+1,\ldots,k+s-1\}$, from (\ref{eq:opt-x-1}) and the convexity of $G_i(\cdot,\xi^l)$ one has
$$
\begin{array}{l}
\langle v_0(x^l,\xi^l),x^{l+1}-x^l\rangle + \frac{1}{2\sigma}\|\lambda^{l+1}\|^2+ \frac{\alpha}{2}\|x^{l+1}-x^l\|^2\\[10pt]
\leq \langle v_0(x^l,\xi^l),\widehat x-x^l\rangle + \frac{1}{2\sigma}\left[ \sum_{i=1}^p[\lambda^l_i+\sigma (G_i(x^l,\xi^l)+ \langle v_i(x^l,\xi^l), \widehat x-x^l \rangle)]_+^2\right]\\[10pt]
\quad\quad + \frac{\alpha}{2}(\|\widehat x-x^l\|^2-\|\widehat x-x^{l+1}\|^2)\\[10pt]
\leq \langle v_0(x^l,\xi^l),\widehat x-x^l\rangle + \frac{1}{2\sigma}\|[\lambda^{l}+\sigma G(\widehat x,\xi^l)]_{+}\|^2+\frac{\alpha}{2}(\|\widehat x-x^l\|^2-\|\widehat x-x^{l+1}\|^2)\\[10pt]
\leq \langle v_0(x^l,\xi^l),\widehat x-x^l\rangle + \frac{1}{2\sigma}\|\lambda^{l}+\sigma G(\widehat x,\xi^l)\|^2+\frac{\alpha}{2}(\|\widehat x-x^l\|^2-\|\widehat x-x^{l+1}\|^2).
\end{array}
$$
Rearranging terms and using Assumption \ref{assu:cons}
we obtain
\[
\begin{array}{ll}
 \frac{1}{2\sigma} \left[\|\lambda^{l+1}\|^2-\|\lambda^l\|^2\right]\\[10pt]
\leq \langle v_0(x^l,\xi^l),\widehat x-x^{l+1}\rangle  + \langle \lambda^l,G(\widehat x,\xi^l)\rangle
+ \frac{\sigma}{2}\|G(\widehat x,\xi^l)\|^2\\[10pt]
\quad\quad + \frac{\alpha}{2}(\|\widehat x-x^l\|^2-\|\widehat x-x^{l+1}\|^2)\\[10pt]
\leq \kappa_fR+ \langle \lambda^l,G(\widehat x,\xi^l)\rangle
+ \frac{\sigma}{2}\nu_g^2+ \frac{\alpha}{2}(\|\widehat x-x^l\|^2-\|\widehat x-x^{l+1}\|^2).
\end{array}
\]
Making  a summation  over $\{k,k+1,\ldots,k+s-1\}$ and taking the conditional expectation, we obtain from
Lemma \ref{lem:aux4} that
\[
\begin{array}{ll}
 \frac{1}{2\sigma} \mathbb{E}\left[\|\lambda^{k+s}\|^2-\|\lambda^k\|^2\,|\, \cF_k\right]\\[10pt]
\leq (\kappa_f R  +  \frac{\sigma}{2}\nu_g^2 )s+ \sum_{l=k}^{k+s-1} \mathbb{E}\left[\langle \lambda^l,G(\widehat x,\xi^l)\rangle\,|\, \cF_k\right]+\frac{\alpha}{2}\|\widehat x-x^k\|^2\\[10pt]
\leq (\kappa_f R  +  \frac{\sigma}{2}\nu_g^2  )s -\varepsilon_0 \sum_{l=0}^{s-1} \mathbb{E}\left[\|\lambda^{k+l}\|\,|\, \cF_k\right]+\frac{\alpha}{2}R^2
\\[10pt]
\leq (\kappa_f R  +  \frac{\sigma}{2}\nu_g^2  )s -\varepsilon_0 \sum_{l=0}^{s-1} \mathbb{E}\left[\|\lambda^{k}\|-\sigma\beta_0l \,|\, \cF_k\right]+\frac{\alpha}{2}R^2
\\[10pt]
 \quad \quad (\mbox{from } \|\lambda^{k+1}\|\geq \|\lambda^k\|-\sigma \beta_0)\\[8pt]
\leq (\kappa_f R  +  \frac{\sigma}{2}\nu_g^2  )s + \varepsilon_0\sigma\beta_0 \frac{s(s-1)}{2}
-\varepsilon_0 s \|\lambda^{k}\| +\frac{\alpha}{2}R^2.
\end{array}
\]
Further, we get from Assumption \ref{assu:cons} and (\ref{eq:theta9}) that
\[
\begin{array}{l}
\mathbb{E}\left[\|\lambda^{k+s}\|^2\,|\, \cF_k\right]\\[10pt]
\leq
\|\lambda^k\|^2+2\sigma(\kappa_f R  +  \frac{\sigma}{2}\nu_g^2 )s
+\varepsilon_0\sigma^2\beta_0s(s-1)-2\varepsilon_0\sigma s \|\lambda^{k}\| +\sigma\alpha R^2\\[10pt]
\leq(\|\lambda^k\|- \frac{\varepsilon_0\sigma}{2}s)^2
+\varepsilon_0\sigma^2 \beta_0s(s-1)+ 2\sigma(\kappa_f R  +  \frac{\sigma}{2}\nu_g^2 )s+\sigma\alpha R^2-\varepsilon_0\sigma s \|\lambda^{k}\| \\[10pt]
\leq(\|\lambda^k\|- \frac{\varepsilon_0\sigma}{2}s)^2
+\varepsilon_0\sigma s[\sigma \beta_0(s-1)+ \frac{2(\kappa_f R  +  \frac{\sigma}{2}\nu_g^2 )}{\varepsilon_0}+\frac{\alpha R^2}{\varepsilon_0s}- \vartheta (\sigma,\alpha,s)]\\[10pt]
\leq (\|\lambda^k\|- \frac{\varepsilon_0\sigma}{2}s)^2.
\end{array}
\]
This, together with Jensen's inequality and the fact that $\|\lambda^k\|\geq  \vartheta (\sigma,\alpha,s)\geq \frac{\varepsilon_0\sigma}{2}s$, implies that
$$
\mathbb{E}\left[\|\lambda^{k+s}\|\,|\, \cF_k\right]\leq
\|\lambda^k\|- \frac{\varepsilon_0\sigma}{2}s.
$$
The proof is completed.
\Halmos\endproof

Let us make some comments on inequality (\ref{eq:7}). This result may seem a bit confusing. From the proof, we actually show that:
the inequality $\mathbb{E}[\|\lambda^{k+s}-\lambda^k\| | \cF_k]\leq s\sigma\beta_0$ holds true under the conditions of Lemma 4; in addition, if $\|\lambda^k\|\geq \vartheta(\sigma,\alpha,s)$, the bound can be improved to $\mathbb{E}[\|\lambda^{k+s}-\lambda^k\| | \cF_k]\leq -s\frac{\sigma\varepsilon_0}{2}$.
However, we state it in the form of (\ref{eq:7}) intentionally.
Since this is only a middle result, our true purpose is to show that the conditions  of the following lemma \citep[Lemma 5]{YMNeely2017} are satisfied for $\|\lambda^k\|$.
\begin{lemma}\label{lem:Yu1}
Let $\{Z_t, t \geq 0\}$ be a discrete time stochastic process adapted to a filtration $\{\cF_t, t\geq
0\}$ with $Z_0 = 0$ and $\cF_0 = \{\emptyset, \Omega\}$. Suppose there exist an integer $t_0 >0$, real constants $\theta>0$, $\delta_{\max}>0$ and $ 0 <\zeta \leq \delta_{\max}$ such that
\[
\begin{array}{rl}
|Z_{t+1}-Z_t| & \leq \delta_{\max},\\[12pt]
\mathbb{E}[Z_{t+t_0}-Z_t\,|\, \cF_t] & \leq \left
\{
\begin{array}{ll}
t_0 \delta_{\max}, & \mbox{if } Z_t < \theta,\\[6pt]
-t_0\zeta, & \mbox{if } Z_t \geq \theta,
\end{array}
\right.
\end{array}
\]
hold for all $t \in \{1,2,\ldots\}.$ Then the following properties are satisfied.
\begin{itemize}
\item[(i)] The  following inequality holds,
\begin{equation}\label{eq:aux1}
\mathbb{E}[Z_t] \leq \theta +t_0 \delta_{\max}+t_0  \frac{4 \delta_{\max}^2}{\zeta}\log \left[  \frac{8 \delta_{\max}^2}{\zeta^2} \right],\ \forall t \in \{1,2,\ldots\}.
\end{equation}
\item[(ii)] For any constant $0 < \mu <1$, we have
    $$
    \Pr\left[Z_t\geq z\right] \leq \mu,\ \forall t \in \{1,2,\ldots\},
    $$
    where
\begin{equation}\label{eq:aux2}
    z=\theta +t_0 \delta_{\max}+t_0  \frac{4 \delta_{\max}^2}{\zeta}\log \left[  \frac{8 \delta_{\max}^2}{\zeta^2}\right]+t_0  \frac{4 \delta_{\max}^2}{\zeta}\log\left( \frac{1}{\mu} \right).
\end{equation}
\end{itemize}
\end{lemma}

It is not difficult to verify that, Lemma \ref{lem:aux5} implies that the conditions of Lemma \ref{lem:Yu1}  are satisfied with respect to $\|\lambda^k\|$
if we take
\[
\theta=\vartheta (\sigma,\alpha,s),\ \delta_{\max}=\sigma \beta_0,\ \zeta = \frac{\sigma}{2}\varepsilon_0,\ t_0=s.
\]
For simplicity,  we define
\[
\psi(\sigma,\alpha,s):= \kappa_0+\kappa_1\frac{\alpha}{s}+\kappa_2\sigma+\kappa_3 \sigma s,\ \phi (\sigma,\alpha,s,\mu):=\psi(\sigma,\alpha,s)+\frac{8\beta_0^2}{\varepsilon_0} \log \left(  \frac{1}{\mu}\right)\sigma s,
\]
where $\kappa_0,\kappa_1,\kappa_2,\kappa_3$ are constants given by
\be\label{eq:kappas}
\kappa_0= \frac{2\kappa_f R}{\varepsilon_0},\,\,
\kappa_1= \frac{ R^2}{\varepsilon_0},\,\,
\kappa_2= \frac{ \nu_g^2}{\varepsilon_0}-\beta_0,\,\,
\kappa_3=\left[2\beta_0 + \frac{\varepsilon_0}{2}+ \frac{8\beta_0^2}{\varepsilon_0}\log  \frac{32\beta_0^2}{\varepsilon_0^2}\right].
\ee
We can also observe that
 $\psi(\sigma,\alpha,s)$ and $\phi (\sigma,\alpha,s,\mu)$ are exactly the same as the right-hand sides of  (\ref{eq:aux1}) and (\ref{eq:aux2}), respectively. Therefore, in view of  Lemma \ref{lem:aux5}, the following lemma is a direct consequence of Lemma  \ref{lem:Yu1}.
\begin{lemma}\label{lem:lambda}
Let Assumptions \ref{assu:compact}--\ref{assu:slater} be satisfied and $s > 0$ be an arbitrary integer. Then, it holds that
\be\label{eq:lambda-Exp}
\Exp[\|\lambda^k\|]\leq \psi(\sigma,\alpha,s).
\ee
Moreover, for any constant $0<\mu<1$, we have
\be\label{eq:lambda-Pr}
\Pr[\|\lambda^k\|\geq \phi(\sigma,\alpha,s,\mu)]\leq\mu.
\ee
\end{lemma}

\section{Expected convergence rates}\label{sec:rates}
In this section, we shall establish the expected convergence rates of SLPMM with respect to constraint violation and objective reduction.

In the following lemma, we derive a bound of the constraints in terms of the averaged iterate $$\hat{x}^K=\frac{1}{K}\sum_{k=0}^{K-1}x^k,$$ where $K$ is a fixed iteration number.
\begin{lemma}\label{lem:cons}
Let  Assumptions \ref{assu:compact}-\ref{assu:moment} be satisfied. Then, if the parameters  satisfy $2\alpha-p\kappa_g^2\sigma>0$, for each $i=1,\ldots,p$ we have
\[
\Exp[g_i(\hat{x}^K)]\leq\frac{1}{\sigma K}\Exp[\lambda^{K}_i]+  \frac{\kappa_g}{\alpha
}(\kappa_f+
\sqrt{p}\nu_g\kappa_g\sigma)+\frac{\sqrt{p}\kappa_g^2}{\alpha K}\sum_{k=0}^{K-1}\Exp[\|\lambda^k\|].
\]
\end{lemma}
\proof{Proof.}
From the definition $\lambda^{k+1}_i=[\lambda^k_i+\sigma ( G_i(x^k,\xi^k)+\langle v_i(x^k,\xi^k), x^{k+1}-x^k\rangle )]_+$, it follows  that
$$
\begin{array}{ll}
\lambda^{k+1}_i
&\geq \lambda^k_i+\sigma( G_i(x^k,\xi^k)+\langle v_i(x^k,\xi^k), x^{k+1}-x^k\rangle )\\[10pt]
& \geq \lambda^k_i+\sigma(G_i(x^k,\xi^k)- \kappa_g\|x^{k+1}-x^k\|).
\end{array}
$$
Using Lemma \ref{lem:aux3}, we have
\be\label{eq:a2}
G_i(x^k,\xi^k)\leq\frac{1}{\sigma}(\lambda^{k+1}_i-\lambda_i^k)+\frac{\kappa_g}{\alpha}(\kappa_f+ \sqrt{p}\kappa_g\|\lambda^k\|+
\sqrt{p}\nu_g\kappa_g\sigma).
\ee
Taking conditional expectation with respect to $\cF_k$, it yields that
\[
g_i(x^k)\leq\frac{1}{\sigma}(\Exp[\lambda^{k+1}_i|\cF_k]-\lambda_i^k)+\frac{\kappa_g}{\alpha}(\kappa_f+ \sqrt{p}\kappa_g\|\lambda^k\|+
\sqrt{p}\nu_g\kappa_g\sigma),
\]
which further gives that
\[
\Exp[g_i(x^k)]\leq\frac{1}{\sigma}(\Exp[\lambda^{k+1}_i]-\Exp[\lambda_i^k])+\frac{\kappa_g}{\alpha}(\kappa_f+ \sqrt{p}\kappa_g\Exp[\|\lambda^k\|]+
\sqrt{p}\nu_g\kappa_g\sigma).
\]
Summing over $\{0,\ldots,K-1\}$ and noticing that $\lambda^0=0$, we obtain
$$
\sum_{k=0}^{K-1} \Exp[g_i(x^k)]\leq
\frac{1}{\sigma}\Exp[\lambda^{K}_i]+  \frac{\kappa_g K}{\alpha
}(\kappa_f+
\sqrt{p}\nu_g\kappa_g\sigma)+\frac{\sqrt{p}\kappa_g^2}{\alpha}\sum_{k=0}^{K-1}\Exp[\|\lambda^k\|].
$$
Therefore, from the convexity of $g_i$ and the definition of $\hat{x}^K$ it follows
\[
\begin{array}{ll}
\Exp[g_i(\hat{x}^K)]&\leq\frac{1}{K}\sum_{k=0}^{K-1} \Exp[g_i(x^k)]\\[10pt]
&\leq\frac{\Exp[\lambda^{K}_i]}{\sigma K}+  \frac{\kappa_g(\kappa_f+
\sqrt{p}\nu_g\kappa_g\sigma)}{\alpha
}+\frac{\sqrt{p}\kappa_g^2}{\alpha K}\sum_{k=0}^{K-1}\Exp[\|\lambda^k\|].
\end{array}
\]
The proof is completed.
\Halmos\endproof

In what follows, we present the bound of the objective reduction in terms of the averaged iterate.
\begin{lemma}\label{lem:obj}
 Let  Assumptions \ref{assu:compact}-\ref{assu:moment} be satisfied. Then, for any $z \in \Phi$,
\[
\Exp[f(\hat{x}^K)]-f(z)\leq\frac{\kappa_f^2}{2\alpha}+\frac{\sigma}{2}\nu_g^2+\frac{\alpha}{2K}R^2.
\]
\end{lemma}
\proof{Proof.}  For any $z \in\Phi$, since $v_0(x^k,\xi^k)\in \partial_x F(x^k,\xi^k)$, we have
$$
\langle v_0(x^k,\xi^k), z-x^k \rangle\leq F(z,\xi^k)-F(x^k,\xi^k).
$$
Then, in view of (\ref{eq:opt-x-1}), one has
\begin{equation}\label{eq:h1}
\begin{array}{ll}
& F(x^k,\xi^k) \\[8pt]
& \leq
 F(z,\xi^k)+\left[\langle v_0(x^k,\xi^k), x^k-x^{k+1}\rangle- \frac{\alpha}{2}\|x^{k+1}-x^k\|^2\right
]\\[6pt]
&\quad +  \frac{1}{2\sigma}\left[\|[\lambda^{k}+\sigma (G(x^k,\xi^k)+V(x^k,\xi^k)(z-x^k))]_+\|^2-\|\lambda^k\|^2\right]\\[12pt]
& \quad -\frac{1}{2\sigma}\left[\|\lambda^{k+1}\|^2-\|\lambda^k\|^2\right] + \frac{\alpha}{2}\left[\|z-x^k\|^2-\|z-x^{k+1}\|^2\right].
\end{array}
\end{equation}
From  Assumption \ref{assu:moment} and the fact that $\langle x,y\rangle\leq \frac{\alpha}{2}\|x\|^2+\frac{1}{2\alpha}\|y\|^2$, we obtain that
\begin{equation}\label{eq:h2}
\langle v_0(x^k,\xi^k),x^k-x^{k+1}\rangle- \frac{\alpha}{2}\|x^{k+1}-x^k\|^2
\leq \frac{1}{2\alpha}\|v_0(x^k,\xi^k)\|^2
\leq \frac{\kappa_f^2}{2\alpha}.
\end{equation}
For every $i=1,\ldots,p$,  we have from $v_i(x^k,\xi^k)\in \partial_x G_i(x^k,\xi^k)$ and $[a]_+^2\leq a^2$ that
$$
[\lambda^k_i+\sigma(G_i(x^k,\xi^k)+\langle v_i(x^k,\xi^k),z-x^k\rangle)]_+^2 \leq [\lambda^k_i+\sigma G_i(z,\xi^k)]^2
$$
and hence
$$
\|[\lambda^k+\sigma(G(x^k,\xi^k)+V(x^k,\xi^k)(z-x^k))]_+\|^2\leq \|\lambda^k+\sigma G(z,\xi^k)\|^2.
$$
Then, we obtain
\begin{equation}\label{eq:h3}
\begin{array}{ll}
\|[\lambda^k+\sigma(G(x^k,\xi^k)+V(x^k,\xi^k)(z-x^k))]_+\|^2-\|\lambda^k\|^2\\[10pt]
\leq 2\sigma \langle \lambda^k,G(z,\xi^k)\rangle+\sigma^2\|G(z,\xi^k)\|^2.
\end{array}
\end{equation}
Substituting (\ref{eq:h2}) and (\ref{eq:h3}) into (\ref{eq:h1}), we get
\be\label{eq:a8}
\begin{array}{ll}
F(x^k,\xi^k) &
\leq F(z,\xi^k)+ \frac{\kappa_f^2}{2\alpha}-\frac{1}{2\sigma}\left[\|\lambda^{k+1}\|^2-\|\lambda^k\|^2\right]
  +   \langle \lambda^k,G(z,\xi^k)\rangle\\[8pt]
&\quad  +\frac{\sigma}{2}\|G(z,\xi^k)\|^2
+ \frac{\alpha}{2}\left[\|z-x^k\|^2-\|z-x^{k+1}\|^2\right].
\end{array}
\ee
Taking conditional expectation with respect to $\cF_k$ and noticing that
\[
\Exp[\langle \lambda^k,G(z,\xi^k)\rangle|\cF_k]=\langle \lambda^k,g(z)\rangle\leq 0,
\]
we have
\[
\begin{array}{ll}
f(x^k)-f(z)&\leq  \frac{\kappa_f^2}{2\alpha}-\frac{1}{2\sigma}\left[\Exp[\|\lambda^{k+1}\|^2|\cF_k]-\|\lambda^k\|^2\right]\\[8pt]
&\quad+\frac{\sigma\nu_g^2}{2}+ \frac{\alpha}{2}\left[\|z-x^k\|^2-\Exp[\|z-x^{k+1}\|^2|\cF_k]\right],
\end{array}
\]
which further gives
\[
\begin{array}{ll}
\Exp[f(x^k)]-f(z)&\leq  \frac{\kappa_f^2}{2\alpha}-\frac{1}{2\sigma}\left[\Exp[\|\lambda^{k+1}\|^2]-\Exp[\|\lambda^k\|^2]\right]\\[8pt]
&\quad+\frac{\sigma\nu_g^2}{2}+ \frac{\alpha}{2}\left[\Exp[\|z-x^k\|^2]-\Exp[\|z-x^{k+1}\|^2]\right].
\end{array}
\]
Making a summation and noticing that $\lambda^0=0$, one has
\[
\sum_{k=0}^{K-1}\Exp[f(x^k)]\leq K\left[f(z)+\frac{\kappa_f^2}{2\alpha}+\frac{\sigma}{2}\nu_g^2\right]+\frac{\alpha}{2}\|z-x^0\|^2.
\]
Therefore, from the convexity of $f$ and the definition of $\hat{x}^K$ it follows
\[
\Exp[f(\hat{x}^K)]\leq\frac{1}{K}\sum_{k=0}^{K-1}\Exp[f(x^k)]\leq f(z)+\frac{\kappa_f^2}{2\alpha}+\frac{\sigma}{2}\nu_g^2+\frac{\alpha}{2K}R^2.
\]
The proof is completed.
\Halmos\endproof

Based on Lemma \ref{lem:cons} and Lemma \ref{lem:obj}, if we take $\alpha=\sqrt{K}$, $\sigma=1/\sqrt{K}$ and $s=\ceil{\sqrt{K}}$,  where $\ceil{a}$ denotes the ceiling function that returns the least integer greater than or equal to $a$, the excepted convergence rates of SLPMM with respect to constraint violation and objective reduction are shown to be
$O(1/\sqrt{K})$ in the following theorem.
\begin{theorem}\label{th:rate}
Let  Assumptions \ref{assu:compact}-\ref{assu:slater} be satisfied. If we take $\alpha=\sqrt{K}$ and $\sigma=1/\sqrt{K}$ in Algorithm \ref{alg:SLPMM}, where $K$ is a fixed iteration number. Then, the following statements hold.
\begin{itemize}
\item[(i)] If  $K>\max\{1,p\kappa_g^2/2\}$, then we have
\[
\Exp[g_i(\hat{x}^K)]\leq\frac{(1+\sqrt{p}\kappa_g^2
)\bar{\kappa}+\kappa_g\kappa_f}{\sqrt{K}}+\frac{(1+\sqrt{p}\kappa_g^2
)\kappa_2+\sqrt{p}\nu_g\kappa_g^2}{K},\quad i=1,\ldots,p,
\]
where $\bar{\kappa}:=\kappa_0+\kappa_1
 +2\kappa_3$ and $\kappa_0, \kappa_1, \kappa_2, \kappa_3$ are defined in (\ref{eq:kappas}).
\item[(ii)] For all $K\geq 1$,
\[
\Exp[f(\hat{x}^K)]-f(x^*)\leq\frac{\kappa_f^2+\nu_g^2+R^2}{2\sqrt{K}},
\]
where $x^*$ is any optimal solution to (\ref{eq:1}).
\end{itemize}
\end{theorem}
\proof{Proof.}
Consider item (i). If $K>p\kappa_g^2/2$, we have $2\alpha-p\kappa_g^2\sigma>0$, then it follows from Lemma \ref{lem:cons}  that
\be\label{eq:a1}
\Exp[g_i(\hat{x}^K)]\leq
\frac{1}{\sigma K}\Exp[\lambda^{K}_i]+  \frac{\kappa_g}{\alpha
}(\kappa_f+
\sqrt{p}\nu_g\kappa_g\sigma)+\frac{\sqrt{p}\kappa_g^2}{\alpha K}\sum_{k=0}^{K-1}\Exp[\|\lambda^k\|].
\ee
If we take $s=\ceil{\sqrt{K}}$,  then from Lemma \ref{lem:lambda} one has
\[
\Exp[\|\lambda^k\|]\leq \psi(\sigma,\alpha,s)=\kappa_0+\kappa_1\frac{\alpha}{s}+\kappa_2
 \sigma+\kappa_3 \sigma s\leq \kappa_0+\kappa_1+\frac{\kappa_2}{\sqrt{K}}+2\kappa_3=\bar{\kappa}+\frac{\kappa_2}{\sqrt{K}}.
\]
Therefore, from $\alpha=\sqrt{K}, \sigma=1/\sqrt{K}$ and (\ref{eq:a1}) we have
\[
\Exp[g_i(\hat{x}^K)]\leq \frac{1}{\sqrt{K}}\left(\bar{\kappa}+\frac{\kappa_2}{\sqrt{K}}\right)+\frac{\kappa_g\kappa_f}{\sqrt{K}}+\frac{\sqrt{p}\nu_g\kappa_g^2}{K}+\frac{\sqrt{p}\kappa_g^2}{\sqrt{K}}\left(\bar{\kappa}+\frac{\kappa_2}{\sqrt{K}}\right),
\]
which verifies item (i).

By taking $z=x^*$ in Lemma \ref{lem:obj}, we derive item (ii)  since
\[
\Exp[f(\hat{x}^K)]-f(x^*)\leq\frac{\kappa_f^2}{2\alpha}+\frac{\sigma}{2}\nu_g^2+\frac{\alpha}{2K}R^2=\frac{\kappa_f^2+\nu_g^2+R^2}{2\sqrt{K}}.
\]
The proof is completed.
\Halmos\endproof

Let us point out that all of the algorithms
 \citep{YMNeely2017,LanZ2016,ABR2021}  have  $O(1/\sqrt{K})$ expected convergence. However, the algorithm \citep{YMNeely2017} is an extension of Zinkevich's online algorithm \citep{Zi2003}, which is a variant of the projection gradient method, and the CSA method \citep{LanZ2016} is  a stochastic counterpart of Polyak's subgradient method \citep{Polyak1967}. When problem (\ref{eq:1}) reduces to a  deterministic problem, these algorithms have at most linear rate of convergence. In contrast, SLPMM becomes the (linearized) proximal method of multipliers, which has an asymptotic superlinear rate of convergence.
Moreover, the iteration complexity analysis \citep{LanZ2016} is based on the selection of stepsizes, which are dependent  on the parameters $R$ , $\kappa_f$ and $\kappa_g$. However, these data are not known beforehand when problem (\ref{eq:1}) is put forward to solve. Note that, in SLPMM the stepsizes $\sigma$ and $\alpha$ are problem-independent.

\section{High probability performance analysis}\label{sec:prob}
In this section, we shall establish the large-deviation properties of SLPMM. By Theorem \ref{th:rate} and  Markov's inequality, we have
for all $\rho_c>0$ and $\rho_o>0$ that
\be\label{eq:d1}
\Pr\left[g_i(\hat{x}^K)\leq \rho_c\left(\frac{(1+\sqrt{p}\kappa_g^2
)\bar{\kappa}+\kappa_g\kappa_f}{\sqrt{K}}+\frac{(1+\sqrt{p}\kappa_g^2
)\kappa_2+\sqrt{p}\nu_g\kappa_g^2}{K}\right)\right]\geq 1-\frac{1}{\rho_c}
\ee
and
\be\label{eq:d2}
\Pr\left[f(\hat{x}^K)-f(x^*)\leq \rho_o\frac{\kappa_f^2+\nu_g^2+R^2}{2\sqrt{K}}\right]\geq 1-\frac{1}{\rho_o}.
\ee
However, these results are very weak.
In the following, we will show that these high probability bounds can be significantly improved.

We introduce the following standard ``light-tail" assumption, see \citep{Lan2016,LanZ2016,LNSY2019} for instance.
\begin{assumption}\label{assu:lt-cons}
There exists a constant $\sigma_c>0$ such that, for any $x\in \cC$,
\[
\Exp[\exp(\|G_i(x,\xi)-g_i(x)\|^2/\sigma_c^2)]\leq\exp(1),\quad i=1,\ldots,p.
\]	
\end{assumption}
From a well-known result \citep[Lemma 4.1]{Lan2020}, under Assumption \ref{assu:lt-cons} one has for any $\rho\geq 0$ and $i=1,\ldots,p$ that
\be\label{eq:a3}
\Pr\left[\frac{1}{K}\sum_{k=0}^{K-1}g_i(x^k)-\frac{1}{K}\sum_{k=0}^{K-1}G_i(x^k,\xi^k)\geq\frac{\rho\sigma_c}{\sqrt{K}}\right]\leq \exp(-\rho^2/3).
\ee

For the sake of readability, we define the following notations,
\[
\theta_1:=\sigma_c+(1+\sqrt{p}\kappa_g^2)\frac{16\beta_0}{\varepsilon_0},\quad \theta_2:=\kappa_g\kappa_f+(1+\sqrt{p}\kappa_g^2)(\kappa_0+\kappa_1+2\kappa_3)
\]
and
\[
\theta_3:=(1+\sqrt{p}\kappa_g^2)\frac{16\beta_0}{\varepsilon_0},\quad \theta_4:=\sqrt{p}\nu_g\kappa_g^2+(1+\sqrt{p}\kappa_g^2)\kappa_2,
\]
in which $\beta_0$ is defined in Lemma \ref{lem:aux5}, $\kappa_0, \kappa_1, \kappa_2, \kappa_3$ are defined in (\ref{eq:kappas}) and other parameters are defined in Assumptions \ref{assu:compact}-\ref{assu:lt-cons}.

We are now read to state  the main result  on constraint violation.
\begin{theorem}\label{th:cons-pr}
Let  Assumptions \ref{assu:compact}-\ref{assu:lt-cons} be satisfied. We take $\alpha=\sqrt{K}$ and $\sigma=1/\sqrt{K}$ in Algorithm \ref{alg:SLPMM}, where $K$ is a fixed iteration number satisfying  $K>\max\{1,p\kappa_g^2/2\}$. Then, for any $\rho\geq 0$ and $i=1,\ldots,p$,
\[
\Pr\left[g_i(\hat{x}^K)\leq \frac{\theta_1\rho+\theta_2+\theta_3\log(K+1)}{\sqrt{K}}+\frac{\theta_4}{K}\right]\geq 1-\exp(-\rho^2/3)-\exp(-\rho).
\]
\end{theorem}
\proof{Proof.}
Summing (\ref{eq:a2}) over $\{0,\ldots,K-1\}$, we have
\[
\frac{1}{K}\sum_{k=0}^{K-1}G_i(x^k,\xi^k)\leq \frac{\lambda_i^K}{\sigma K}+\frac{\kappa_g(\kappa_f+\sqrt{p}\nu_g\kappa_g\sigma)}{\alpha}+\frac{\sqrt{p}\kappa_g^2}{\alpha K}\sum_{k=0}^{K-1}\|\lambda^k\|.
\]
Noticing that $\alpha=\sqrt{K}$, $\sigma=1/\sqrt{K}$ and $g_i(\hat{x}^K)\leq \frac{1}{K}\sum_{k=0}^{K-1}g_i(x^k)$, one has
\be\label{eq:a4}
g_i(\hat{x}^K)\leq  \frac{1}{K}\sum_{k=0}^{K-1}[g_i(x^k)-G_i(x^k,\xi^k)]+\frac{\lambda_i^K}{\sqrt{K}}+\frac{\kappa_g\kappa_f}{\sqrt{K}}+\frac{\sqrt{p}\nu_g\kappa_g^2}{K}+\frac{\sqrt{p}\kappa_g^2}{K^{3/2}}\sum_{k=0}^{K-1}\|\lambda^k\|.
\ee
We next consider the probability bound of $\lambda^k$. From (\ref{eq:lambda-Pr}), it follows that
\[
\Pr[\|\lambda^k\|\geq \phi(\sigma,\alpha,s,\mu)]\leq\mu,\quad k=0,1\ldots,K.
\]
If we take $s=\ceil{\sqrt{K}}$ and $\mu=\exp(-\rho)/(K+1)$, then
\[
\begin{array}{ll}
\phi(\sigma,\alpha,s,\mu)&= \kappa_0+\kappa_1\frac{\alpha}{s}+\kappa_2\sigma+\kappa_3 \sigma s+\frac{8\beta_0^2}{\varepsilon_0} \log \left(  \frac{1}{\mu}\right)\sigma s\\[10pt]
&\leq \kappa_0+\kappa_1+\frac{\kappa_2}{\sqrt{K}}+2\kappa_3+\frac{16\beta_0^2}{\varepsilon_0}(\rho+\log(K+1))
\end{array}
\]
and hence for all $k=0,1,\ldots,K$,
\be\label{eq:a5}
\Pr[\|\lambda^k\|\geq \kappa_0+\kappa_1+\frac{\kappa_2}{\sqrt{K}}+2\kappa_3+\frac{16\beta_0^2}{\varepsilon_0}(\rho+\log(K+1))]\leq\frac{\exp(-\rho)}{K+1}.
\ee
Using (\ref{eq:a3}) and (\ref{eq:a5}) in (\ref{eq:a4}), we conclude that
\[
\begin{array}{ll}
\Pr\left[ g_i(\hat{x}^K)\geq \frac{\rho(\sigma_c+(1+\sqrt{p}\kappa_g^2)\frac{16\beta_0}{\varepsilon_0})}{\sqrt{K}}+\frac{\kappa_g\kappa_f+(1+\sqrt{p}\kappa_g^2)(\kappa_0+\kappa_1+2\kappa_3)}{\sqrt{K}}\right.\\[15pt]
\quad\quad\left.+\frac{(1+\sqrt{p}\kappa_g^2)\frac{16\beta_0}{\varepsilon_0}\log(K+1)}{\sqrt{K}}+\frac{\sqrt{p}\nu_g\kappa_g^2+(1+\sqrt{p}\kappa_g^2)\kappa_2}{K}
\right]\leq \exp(-\rho^2/3)+\exp(-\rho).
\end{array}
\]
The proof is completed.
\Halmos\endproof

In view of Theorem \ref{th:cons-pr}, if we take $\rho=\log(K)$, then we have
\[
\Pr\left[g_i(\hat{x}^K)\leq O\left(\frac{\log(K)}{\sqrt{K}}\right)\right]\geq 1-\frac{1}{K^{2/3}}-\frac{1}{K}.
\]

We next make the following  ``light-tail" assumption with respect to the objective function.
\begin{assumption}\label{assu:lt-obj}
There exists a constant $\sigma_o>0$ such that, for any $x\in \cC$,
\[
\Exp[\exp(\|F(x,\xi)-f(x)\|^2/\sigma_o^2)]\leq\exp(1).
\]	
\end{assumption}
Similar to (\ref{eq:a3}), under Assumption \ref{assu:lt-obj} one has for any $\rho\geq 0$ that
\be\label{eq:a6}
\Pr\left[\frac{1}{K}\sum_{k=0}^{K-1}f(x^k)-\frac{1}{K}\sum_{k=0}^{K-1}F(x^k,\xi^k)\geq\frac{\rho\sigma_o}{\sqrt{K}}\right]\leq \exp(-\rho^2/3)
\ee
and
\be\label{eq:a7}
\Pr\left[\frac{1}{K}\sum_{k=0}^{K-1}F(z,\xi^k)-\frac{1}{K}\sum_{k=0}^{K-1}f(z)\geq\frac{\rho\sigma_o}{\sqrt{K}}\right]\leq \exp(-\rho^2/3)
\ee
for all $z\in \cC$.

The following  lemma is from \citep[Lemma 9]{YMNeely2017}.
\begin{lemma}\label{lem:Yu2}
Let $\{Z_t,t\geq 0\}$ be a supermartingale adapted to a filtration $\{\cF_t,t\geq 0\}$ with $Z_0=0$ and $\cF_0=\{\emptyset, \Omega\}$, i.e. $\mathbb E[Z_{t+1}\,|\, \cF_t]\leq Z_t$, $\forall t \geq 0$. Suppose there exists a constant $c>0$ such that $\{|Z_{t+1}-Z_t|>c\}\subseteq \{Y_t>0\}$,  $\forall t\geq 0$, where each $Y_t$ is  adapted to $\cF_t$.  Then, for all $z>0$, we have
$$
\Pr[Z_t\geq z] \leq e^{-z^2/(2tc^2)}+ \sum_{j=0}^{t-1} \Pr[Y_j>0],\ \forall t \geq 1.
$$
\end{lemma}

For any fixed $z \in \Phi$, by taking $Z_t:= \sum_{k=0}^{t-1} \langle \lambda^k, G(z, \xi^k) \rangle$ in Lemma \ref{lem:Yu2} we obtain  the following lemma.
\begin{lemma}\label{lem:l13}
For any fixed $z \in \Phi$ and an arbitrary constant $c>0$, let $Z_0:=0$ and $Z_t:= \sum_{k=0}^{t-1} \langle \lambda^k, G(z, \xi^k) \rangle$ for $t\geq 1$. Let  $\cF_0=\{\emptyset, \Omega\}$   and $Y_t:=\|\lambda^{t}\|-c/\nu_g$ for all $t\geq 0$. Then, for all $\gamma>0$, we have
$$
\Pr[Z_t\geq \gamma] \leq e^{-\gamma^2/(2tc^2)}+ \sum_{j=0}^{t-1} \Pr[Y_j>0],\ \forall t \geq 1.
$$
\end{lemma}
\proof{Proof.}  It is simple to check that $\{Z_t\}$ and $\{Y_t\}$ are both adapted to $\{\cF_t, t\geq 0\}$. Now we prove that
$\{Z_t\}$ is a supermartingale. Since
$
Z_{t+1}=Z_t+\langle \lambda^{t}, G(z, \xi^{t})\rangle,
$
we have
$$
\begin{array}{ll}
\mathbb{E}[Z_{t+1}\,|\, \cF_t]&=\mathbb{E}[ Z_t+\langle \lambda^{t}, G(z, \xi^{t})\rangle\,|\, \cF_t]\\[8pt]
& =Z_t+\langle \lambda^{t}, \mathbb{E}[G(z, \xi^{t})\,|\, \cF_t]\rangle\\[8pt]
&=Z_t+\langle \lambda^{t}, g(z)\rangle\\[8pt]
&\leq Z_t,
\end{array}
$$
which follows from $\lambda^t\in\cF_t$, $\lambda^t\geq 0$ and $g(z)\leq 0$. Thus, we obtain that $\{Z_t\}$ is a supermartingale.

From Assumption \ref{assu:cons}, we get
$$
|Z_{t+1}-Z_t|=|\langle \lambda^{t}, G(z, \xi^{t} \rangle| \leq \nu_g\|\lambda^{t}\|.
$$
This implies that $\|\lambda^{t}\|> c/\nu_g$ if $|Z_{t+1}-Z_t| >c$ and hence
$$
\{|Z_{t+1}-Z_t| >c\} \subseteq \{Y_t>0\}.
$$
Therefore, we can observe that the conditions of Lemma \ref{lem:Yu2} are satisfied, and hence the claim is obtained.
\Halmos\endproof

Finally, we establish a high probability objective reduction bound in the following theorem.
\begin{theorem}\label{th:obj-pr}
Let  Assumptions \ref{assu:compact}-\ref{assu:slater} and \ref{assu:lt-obj} be satisfied. We take $\alpha=\sqrt{K}$ and $\sigma=1/\sqrt{K}$ in Algorithm \ref{alg:SLPMM}, where $K\geq 1$ is a fixed iteration number. Then, for any $\rho\geq 0$,
\[
\begin{array}{ll}
\displaystyle\Pr\left[f(\hat{x}^K)-f(x^*)\leq\sqrt{2\rho}\nu_g\left(\frac{\kappa_0+\kappa_1+2\kappa_3}{\sqrt{K}}+\frac{\frac{16\beta_0^2}{\varepsilon_0}(\rho+\log(K))}{\sqrt{K}}+\frac{\kappa_2}{K}\right)\right.\\[15pt]
\quad\quad \displaystyle\left.+ \frac{2\sigma_0\rho}{\sqrt{K}}+\frac{\theta_5}{\sqrt{K}} \right]\geq 1-2\exp(-\rho^2/3)-2\exp(-\rho),
\end{array}
\]
where $x^*$ is any fixed optimal solution to (\ref{eq:1}),
$
\theta_5:=(\kappa_f^2+\nu_g^2+R^2)/2$, $\beta_0$ is defined in Lemma \ref{lem:aux5} and $\kappa_0, \kappa_1, \kappa_2, \kappa_3$ are defined in (\ref{eq:kappas}).

\end{theorem}
\proof{Proof.}
For any $z\in\Phi$, summing (\ref{eq:a8}) over $\{0,\ldots,K-1\}$ and using the facts that $\lambda^0=0$, $\|G(z,\xi^k)\|^2\leq \nu_g^2$ and $\|z-x^0\|^2\leq R^2$, we have
\[
\frac{1}{K}\sum_{k=0}^{K-1}F(x^k,\xi^k)\leq \frac{1}{K}\sum_{k=0}^{K-1}F(z,\xi^k)+\frac{1}{K}\sum_{k=0}^{K-1}\langle \lambda^k,
G(z, \xi^k)\rangle +\frac{\kappa_f^2+\nu_g^2+R^2}{2\sqrt{K}}.
\]
Then, it follows from $f(\hat{x}^K)\leq \frac{1}{K}\sum_{k=0}^{K-1}f(x^k)$ that
\be\label{eq:b1}
\begin{array}{ll}
&f(\hat{x}^K)-f(z)\\[10pt]
&\leq \displaystyle\frac{1}{K}\sum_{k=0}^{K-1}[f(x^k)-F(x^k,\xi^k)]+\frac{1}{K}\sum_{k=0}^{K-1}[F(z,\xi^k)-f(z)]\\[15pt]
&\quad \displaystyle+\frac{1}{K}\sum_{k=0}^{K-1}\langle \lambda^k,
G(z, \xi^k)\rangle+\frac{\kappa_f^2+\nu_g^2+R^2}{2\sqrt{K}}.
\end{array}
\ee
By Lemma \ref{lem:l13}, for any $c>0$ and $\gamma>0$ we have
\[
\Pr\left[\frac{1}{K}\sum_{k=0}^{K-1}\langle \lambda^k,
G(z, \xi^k)\rangle\geq \frac{\gamma}{K}\right] \leq \exp(-\gamma^2/(2Kc^2))+ \sum_{k=0}^{K-1}\Pr[\|\lambda^k\|\geq c/\nu_g].
\]
Let us take $s=\ceil{\sqrt{K}}$ and $\mu=\exp(-\rho)/K$, then
\[
\phi(\sigma,\alpha,s,\mu)\leq \kappa_0+\kappa_1+\frac{\kappa_2}{\sqrt{K}}+2\kappa_3+\frac{16\beta_0^2}{\varepsilon_0}(\rho+\log(K)).
\]
If we take $c=\nu_g\phi(\sigma,\alpha,s,\mu)$, then from (\ref{eq:lambda-Pr}) we obtain
\[
\sum_{k=0}^{K-1}\Pr[\|\lambda^k\|\geq c/\nu_g]\leq K\mu=\exp(-\rho).
\]
Moreover, let us take $\gamma=\sqrt{2\rho K} c$, then
\[
\frac{\gamma}{K}
=\sqrt{2\rho}\nu_g\left(\frac{\kappa_0+\kappa_1+2\kappa_3}{\sqrt{K}}+\frac{\frac{16\beta_0^2}{\varepsilon_0}(\rho+\log(K))}{\sqrt{K}}+\frac{\kappa_2}{K}\right)
\]
and hence
\be\label{eq:b2}
\begin{array}{ll}
\Pr\left[\frac{1}{K}\sum_{k=0}^{K-1}\langle \lambda^k,
G(z, \xi^k)\rangle\geq \sqrt{2\rho}\nu_g\left(\frac{\kappa_0+\kappa_1+2\kappa_3}{\sqrt{K}}+\frac{\frac{16\beta_0^2}{\varepsilon_0}(\rho+\log(K))}{\sqrt{K}}+\frac{\kappa_2}{K}\right)\right] \\[15pt]
\leq 2\exp(-\rho).
\end{array}
\ee
Using (\ref{eq:a6}), (\ref{eq:a7}) and (\ref{eq:b2}) in (\ref{eq:b1}), one has
\[
\begin{array}{ll}
\Pr\left[f(\hat{x}^K)-f(z)\geq\frac{2\sigma_0\rho}{\sqrt{K}}+\sqrt{2\rho}\nu_g\left(\frac{\kappa_0+\kappa_1+2\kappa_3}{\sqrt{K}}+\frac{\frac{16\beta_0^2}{\varepsilon_0}(\rho+\log(K))}{\sqrt{K}}+\frac{\kappa_2}{K}\right)\right.\\[15pt]
\quad\quad\left.+\frac{\kappa_f^2+\nu_g^2+R^2}{2\sqrt{K}} \right]\leq 2\exp(-\rho^2/3)+2\exp(-\rho).
\end{array}
\]
The claim is derived by taking $z=x^*$ in the above inequality.
\Halmos\endproof
In view of Theorem \ref{th:obj-pr}, if we take $\rho=\log(K)$, then we have
\[
\Pr\left[f(\hat{x}^K)-f(x^*)\leq O\left(\frac{\log^{3/2}(K)}{\sqrt{K}}\right)\right]\geq 1-\frac{2}{K^{2/3}}-\frac{2}{K}.
\]
In contrast to (\ref{eq:d1}) and (\ref{eq:d2}), we can observe that the results in Theorem \ref{th:cons-pr} and \ref{th:obj-pr} are much finer.

\section{Preliminary numerical experiments}\label{sec:num}
In this section, we demonstrate the efficiency of the proposed stochastic linearized proximal method of multipliers on two preliminary numerical problems. All numerical experiments are carried out using MATLAB R2020a on a desktop computer with Intel(R) Xeon(R) E-2124G 3.40GHz and 32GB memory. The MATLAB code and test problems can be found on
\url{https://bitbucket.org/Xiantao_Xiao/SLPMM}. All reported time is wall-clock time in seconds.
\subsection{Solving subproblems}\label{sec:subp}
This subsection focuses on solving the  subproblem (\ref{xna}) in SLPMM, that is
\[
\begin{array}{l}
 x^{k+1}= \argmin\limits_{x \in \cC} \,\left\{ \cL^k_{\sigma }(x,\lambda^k) +\frac{\alpha}{2}\|x-x^k\|^2\right\}.
\end{array}
\]
This problem is equivalent to
\begin{equation}\label{eq:general-subp}
 \min_{x \in \cC} \phi(x):=\frac{1}{2}\sum_{i=1}^p[a_i^Tx+b_i]_+^2+\frac{1}{2}\|x\|^2+c^Tx,
\end{equation}
where
\[
a_i:=\sqrt{\frac{\sigma}{\alpha}}v_i(x^k,\xi^k),\ b_i:=\frac{\lambda_i}{\sqrt{\sigma\alpha}}+\sqrt{\frac{\sigma}{\alpha}}G_i(x^k,\xi^k)-\left\langle\sqrt{\frac{\sigma}{\alpha}}v_i(x^k,\xi^k),x^k\right\rangle
\]
and $c:=v_0(x^k,\xi^k)/\alpha-x^k$.
Since $\phi$ is obviously strongly convex,  we could apply the following popular Nesterov's accelerated gradient method to solve (\ref{eq:general-subp}).
\mbox{}\\[4pt]
{\bf APG}: Nesterov's accelerated projected gradient method  for   (\ref{eq:general-subp}).
\begin{description}
\item[Step 0 ] Input $x^0\in \cC$  and $\eta>1$. Set $y^0=x^0$, $L_{-1}=1$ and $t:=0$.
\item[Step 1]
Set
\[
x^{t+1}=T_{L_t}(y^t),
\]
where $T_{L}(y):=\Pi_{\cC}[y-\frac{1}{L}\nabla \phi(y)]$, the stepsize $L_t=L_{t-1}\eta^{i_t}$ and $i_t$ is the smallest nonnegative integer satisfies the following condition
\[
\begin{array}{ll}
\phi(T_{L_{t-1}\eta^{i_t}}(y^t))&\leq \phi(y^t)+\langle\nabla \phi(y^t),T_{L_{t-1}\eta^{i_t}}(y^t)-y^t\rangle\\[10pt]
&\quad\quad+\frac{L_{t-1}\eta^{i_t}}{2}\|T_{L_{t-1}\eta^{i_t}}(y^t)-y^t\|^2.
\end{array}
\]
\item[ Step 2] Compute
\[
y^{t+1}=x^{t+1}+\frac{t}{t+3}(x^{t+1}-x^t).
\]
\item[ Step 3] Set $t:=t+1$  and go to Step 1.
\end{description}

A well-known convergence result of the above method is that, if $\phi$ is $\mu$-strongly convex and $\nabla \phi$ is $L$-Lipschitz continuous, then
$\phi(x^t)-\phi(x^*)\leq O\left((1-\sqrt{\mu/L})^t\right)$. See \citep{Beck2017} for a detailed discussion on this topic. Here, we assume that the set $\cC$ is simple such  that the projection $\Pi_{\cC}$ can be efficiently computed. For example, if
\[\cC:=\left\{x\in\R^n:\sum_{i=1}^nx_i=1,\ x\geq 0\right\},\]  the projection $\Pi_{\cC}$ can be  computed by the method proposed in \citep{WL2015}.

When  $\cC$ is $\R^n$ or a polyhedron, the subproblem is equivalent to a convex quadratic programming (QP) problem as
\[
\begin{array}{ll}
	\min\limits_{x,y}\quad & \displaystyle\frac{1}{2}\sum_{i=1}^py_i^2+\frac{1}{2}\|x\|^2+c^Tx\\[8pt]
	\textrm{s.t.}\quad & a_i^Tx+b_i-y_i\leq 0,\ i=1,2,\ldots,p,\\[5pt]
	& x\in\cC,\quad y\geq 0.
\end{array}
\]
In this case, the subproblem  can also be solved by a QP solver.
Let us also mention that, if $p=1$,  the closed form of the stationary point to the objective function in Problem (\ref{eq:general-subp}) is  given by
\[
\tilde{x}=\left\{
\begin{array}{ll}
-c,\quad&\mbox{if}\ -a_1^Tc+b_1\leq 0,\\[8pt]
-(b_1a_1+c)+\frac{a_1^T(b_1a_1+c)a_1}{1+a_1^Ta_1},\quad&\mbox{otherwise}.
\end{array}
\right.
\]
Then, $\tilde{x}$ is the unique optimal solution if it lies in the interior of $\cC$.

\subsection{Neyman-Pearson classification}
For a classifier $h$ to predict $1$ and $-1$, let us define the type I error (misclassifying class -1 as 1) and type II error (misclassifying class 1 as -1) respectively by
\[
\mbox{type I error}:=\mathbb{E}[\varphi(-bh(a))|b=-1],\quad\mbox{type II error}:=\mathbb{E}[\varphi(-bh(a))|b=1],
\]
where $\varphi$ is some merit function. Unlike the conventional binary classification in machine learning, the Neyman-Pearson (NP) classification paradigm is developed to learn a classifier by minimizing type II error with type I error being below a user-specified level $\tau>0$,
see  \citep{TFZ2016} and  references therein.
In specific, for a given class $\mathcal{H}$ of classifiers,  the NP classification is to solve the following problem
\[
\begin{array}{ll}
\min\limits_{h\in\mathcal{H}}&\mathbb{E}[\varphi(-bh(a))|b=1]\\[5pt]
\mbox{s.t.}\quad &\mathbb{E}[\varphi(-bh(a))|b=-1]\leq \tau.
\end{array}
\]
In what follows, we consider its empirical risk minimization counterpart. Suppose that a labeled training dataset $\{a_i\}_{i=1}^N$ consists of the positive set $\{a^0_i\}_{i=1}^{N_0}$ and the negative set $\{a^1_i\}_{i=1}^{N_1}$. The associated empirical NP classification problem is
\begin{equation}\label{prob:NP}
\begin{array}{ll}
\min\limits_x&f(x):=\frac{1}{N_0}\sum_{i=1}^{N_0}\ell(x^Ta_i^0)\\[8pt]
\mbox{s.t.}\quad &g(x):=\frac{1}{N_1}\sum_{i=1}^{N_1}\ell(-x^Ta_i^1)-\tau\leq 0,
\end{array}
\end{equation}
where $\ell(\cdot)$ is a loss function, e.g., logistic loss $\ell(y):=\log(1+\exp(-y))$.

The datasets tested in our numerical comparison are summarized in Table \ref{table:datasets}. The datasets for multi-class classification have been manually divided into two types. For example, the MNIST dataset is used for classifying odd and even digits.
\begin{table}[!htp]
\centering
\caption{Datasets used in  Neyman-Pearson classification}
\begin{tabular}{|c||c|c|c|c|}
\hline
Dataset & Data  $N$ & Variable $n$ & {Density} & Reference \\
\hline
$\mathtt{gisette}$ & 6000 & 5000 & 12.97\% &\citep{gisette}\\[2pt]
$\mathtt{CINA}$ & 16033 & 132 & 29.56\% &\citep{CINA}\\[2pt]
$\mathtt{MNIST}$ & 60000 & 784 & 19.12\% &\citep{MNIST}\\
\hline
\end{tabular}
\label{table:datasets}
\end{table}



In the following experiment, we show the performance of SLPMM compared with CSA \citep{LanZ2016},  PSG \citep{Xiao2019}, YNW \citep{YMNeely2017} and APriD \citep{YX2022}. For all five methods, we use an efficient mini-batch strategy, that is, at each iteration the stochastic gradients of the objective function and the constraint function are computed, respectively, by
\[
v_0^k:=\frac{1}{|\cN_0^k|}\sum_{i\in\cN_0^k}\nabla f_i(x^k),\quad v_1^k:=\frac{1}{|\cN_1^k|}\sum_{i\in\cN_1^k}\nabla g_i(x^k),
\]
where $f_i(x):=\ell(x^Ta_i^0), i=1,\ldots,N_0$ and $g_i(x):=\ell(-x^Ta_i^1), i=1,\ldots,N_1$. Here, the sets $\cN_0^k$ and $\cN_1^k$ are randomly chosen from the index sets $\{1,\ldots,N_0\}$ and $\{1,\ldots,N_1\}$, respectively. The batch sizes $|\cN_0^k|$ and  $|\cN_1^k|$ are fixed to $1\%$ of the data sizes $N_0$ and $N_1$, respectively.
We choose $x^0=0$ as the initial point. The parameter $\tau$ is set to 1. The parameters in SLPMM is chosen as $\alpha=\sqrt{K}$ and $\sigma=1/\sqrt{K}$. The maximum number of iterations is set to $K=3000$.

In Figure \ref{figure:gisette}, Figure \ref{figure:CINA} and Figure \ref{figure:MNIST}, we show the performance of all methods for solving the empirical NP classification problem with logistic loss. In each figure, the pictures (a) and (b) show the changes of the objective value and the constraint value with respect to \textit{epochs}, and the pictures (c) and (d) represent the changes of the objective value and the constraint value with respect to \textit{cputime}. Here, in (a) and  (c) the horizontal dashed line represents a reference optimal objective value which is computed by the built-in MATLAB function \texttt{fmincon}. Moreover, one epoch denotes a full pass over a dataset. The results are averaged over 10 independent runs.

Generally, we can observe that the behaviors of CSA, PSG and YNW are similar since all of them are stochastic first-order methods.
 SLPMM obviously outperforms these three methods by combining the evaluations of both objective decreasing and constraint violation. In particular, the results demonstrate that SLPMM converges obviously faster than CSA and PSG both with respect to epochs and cputime.    Our results also show that PSG usually generates solutions which are failed to satisfy the constraint. In contrast, CSA always gives feasible solutions, but the objective values are far from optimal. Finally, the performance of APriD is very different from the others. The total performance of APriD seems better than the others. However, the curves of APriD oscillate heavily even for the average of 10 runs, and the issue is  much worse for each independent run.

\begin{figure}[htp]
\centering
\setlength{\belowcaptionskip}{-6pt}
\begin{tabular}{cccc}
\subfloat[$\mathtt{objective \slash epochs}$]{
\includegraphics[width=5.5cm]{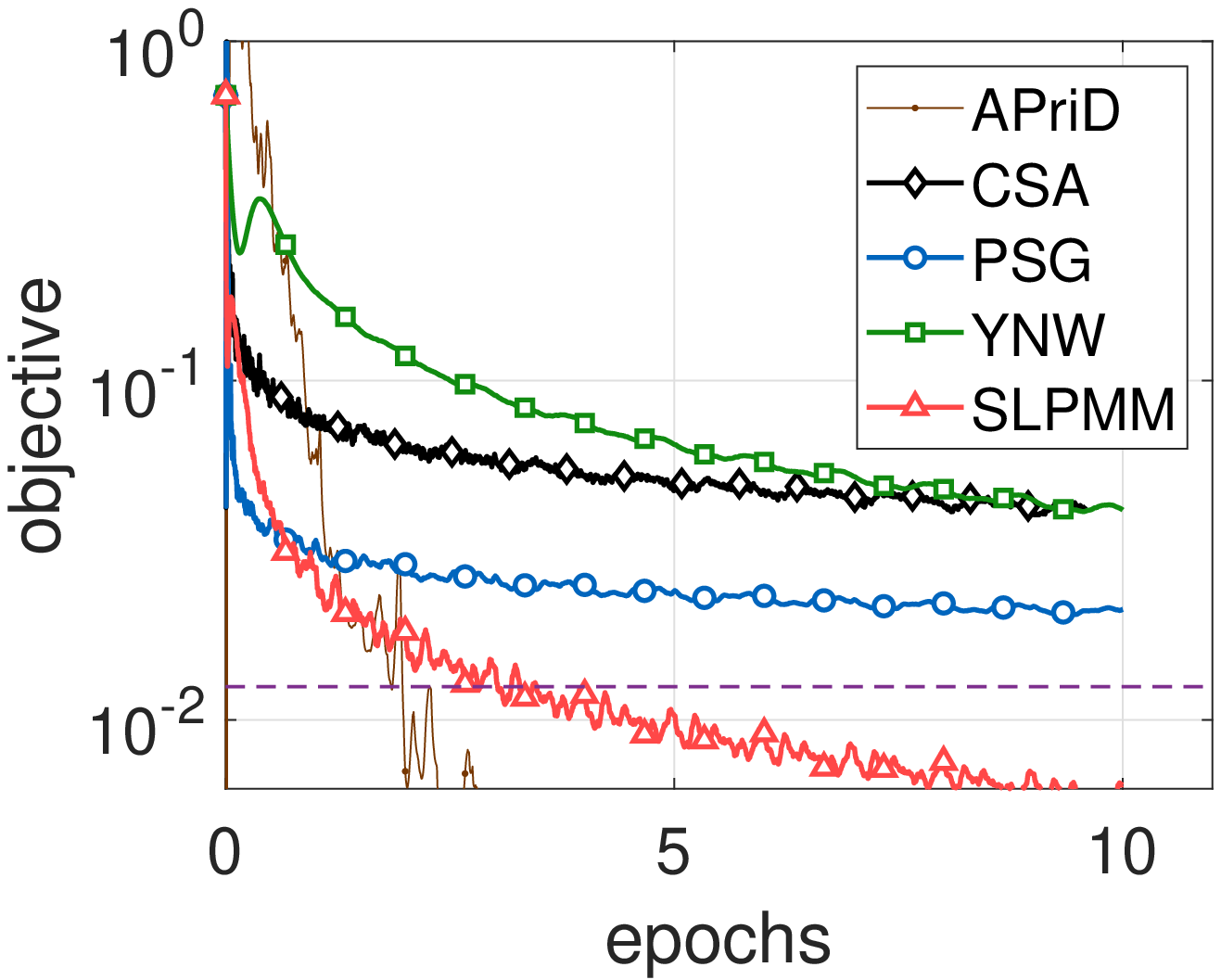}} &
\subfloat[$\mathtt{constraint \slash epochs}$]{
\includegraphics[width=5.5cm]{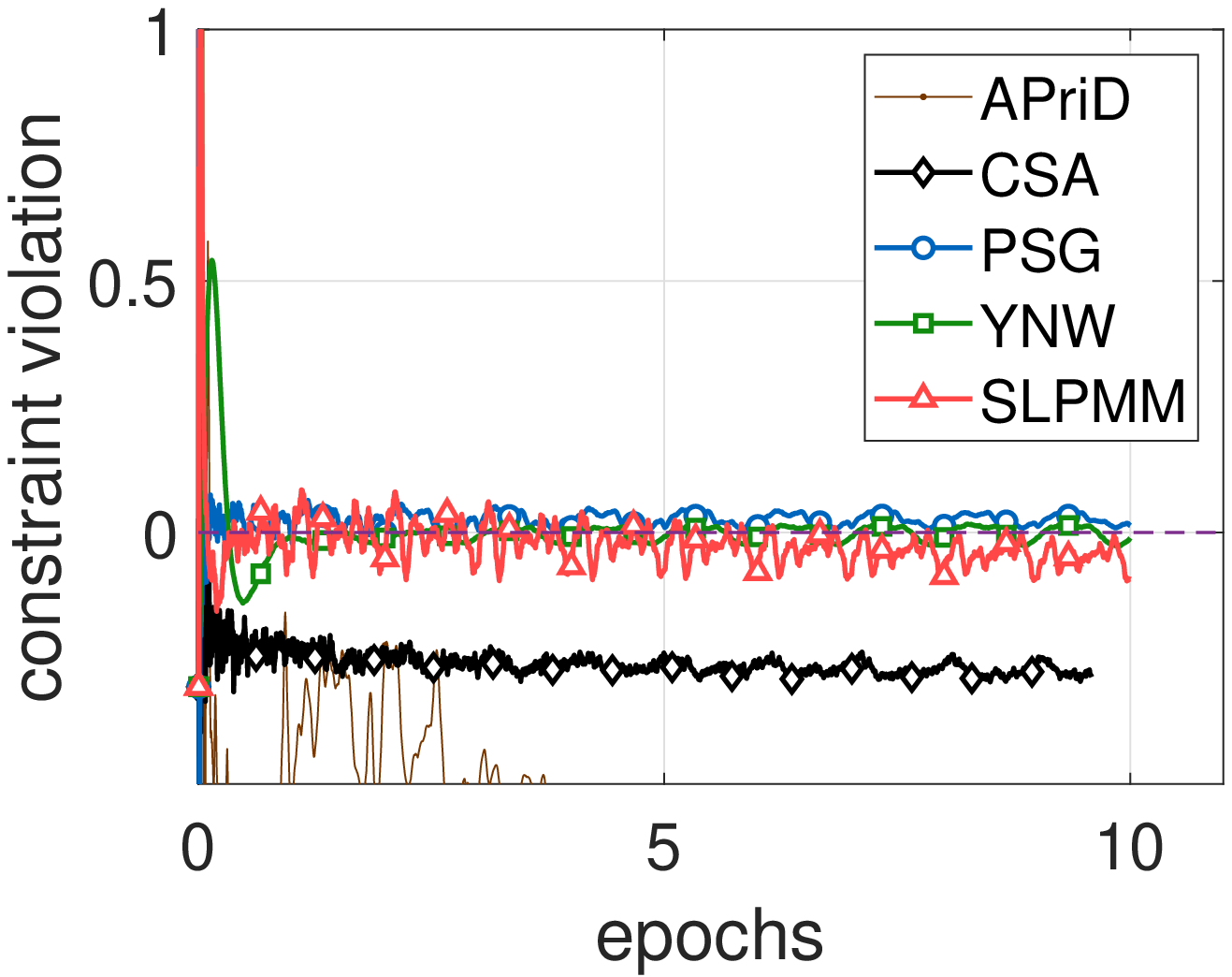}} &\\
\subfloat[$\mathtt{objective \slash cputime}$]{
\includegraphics[width=5.5cm]{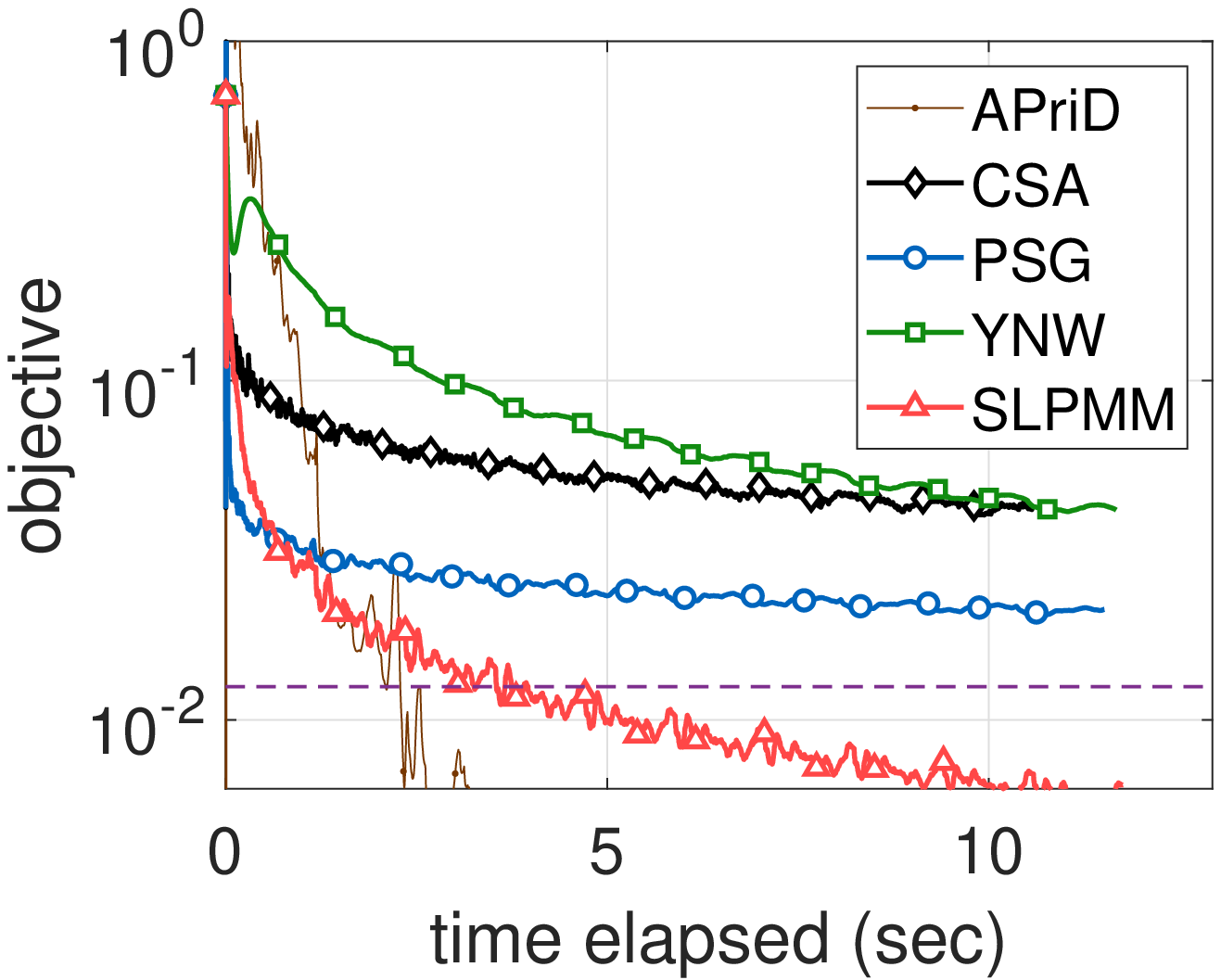}} &
\subfloat[$\mathtt{constraint \slash cputime}$]{
\includegraphics[width=5.5cm]{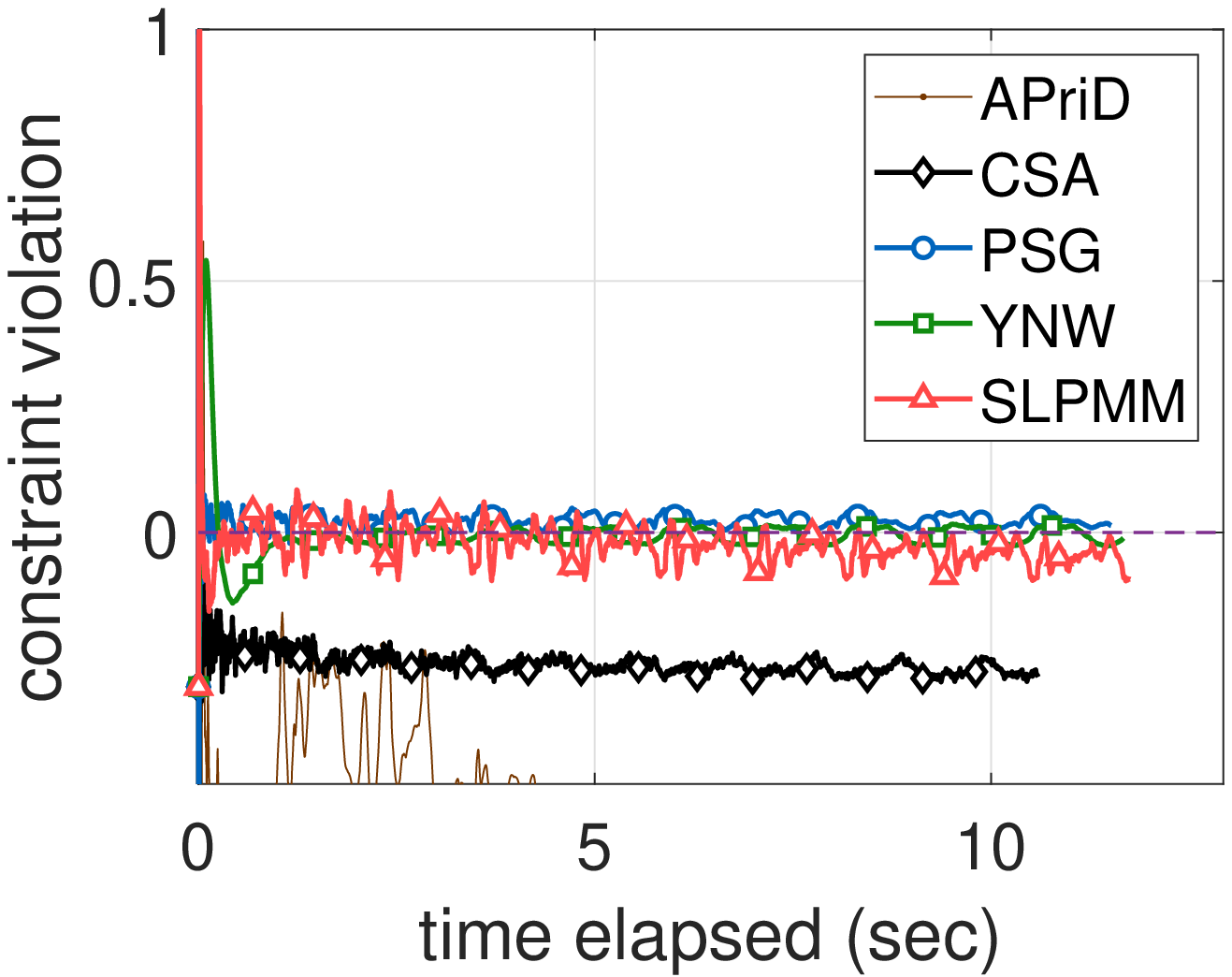}}
\end{tabular}
\caption{Comparison of algorithms on $\mathtt{gisette}$ for Neyman-Pearson classification.}
\label{figure:gisette}
\end{figure}
\begin{figure}[htp]
\centering
\setlength{\belowcaptionskip}{-6pt}
\begin{tabular}{cccc}
\subfloat[$\mathtt{objective \slash epochs}$]{
\includegraphics[width=5.5cm]{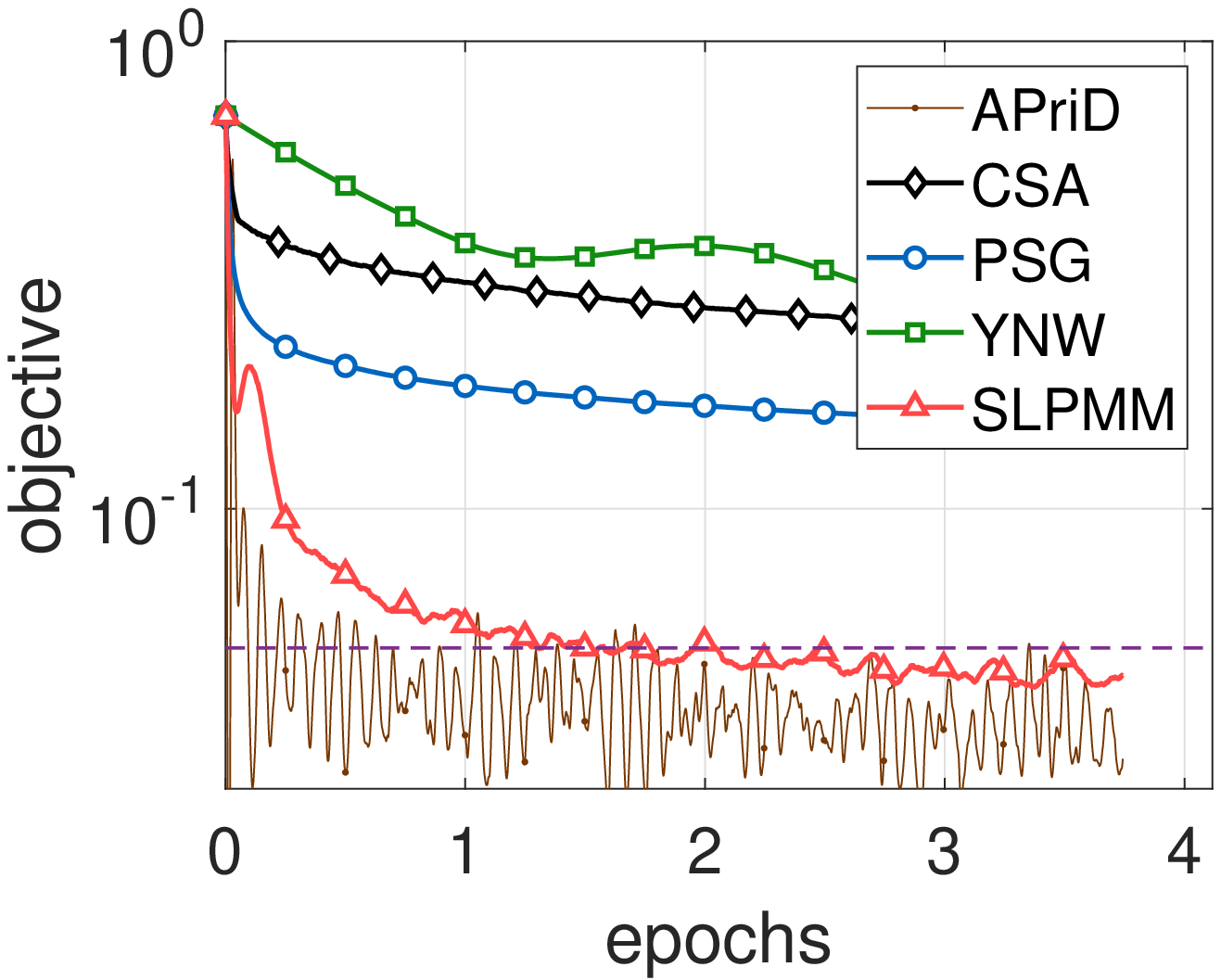}} &
\subfloat[$\mathtt{constraint \slash epochs}$]{
\includegraphics[width=5.5cm]{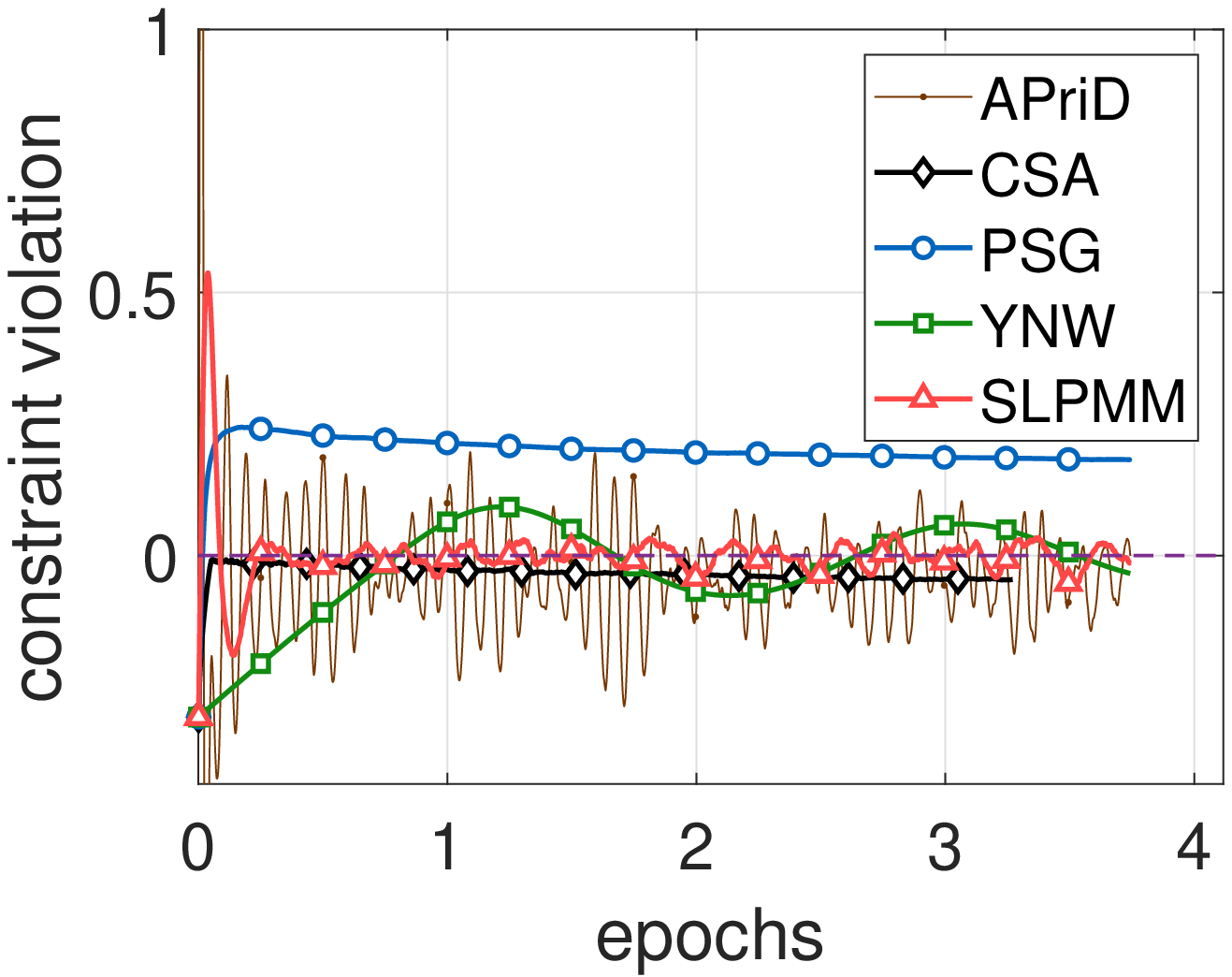}} &\\
\subfloat[$\mathtt{objective \slash cputime}$]{
\includegraphics[width=5.5cm]{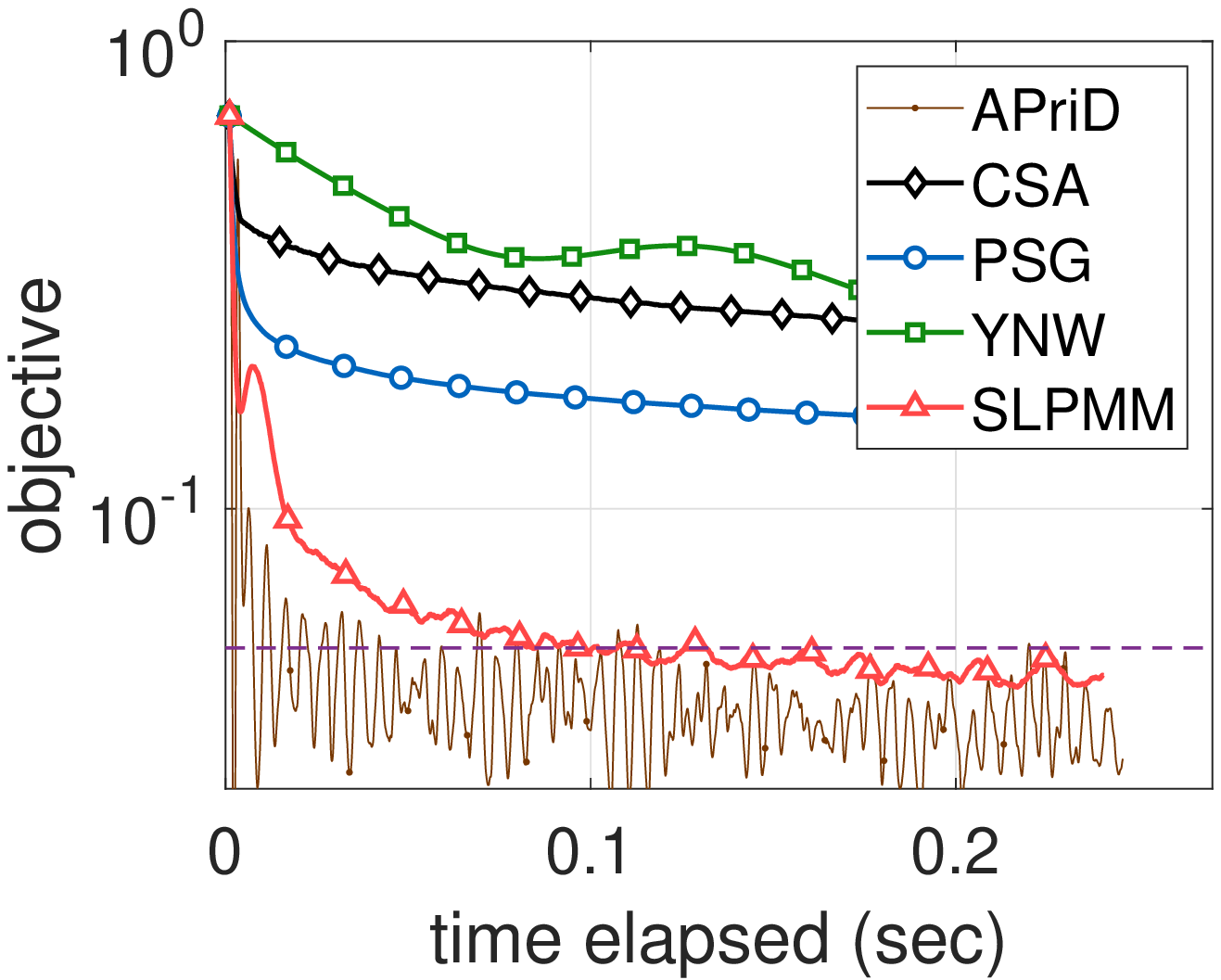}} &
\subfloat[$\mathtt{constraint \slash cputime}$]{
\includegraphics[width=5.5cm]{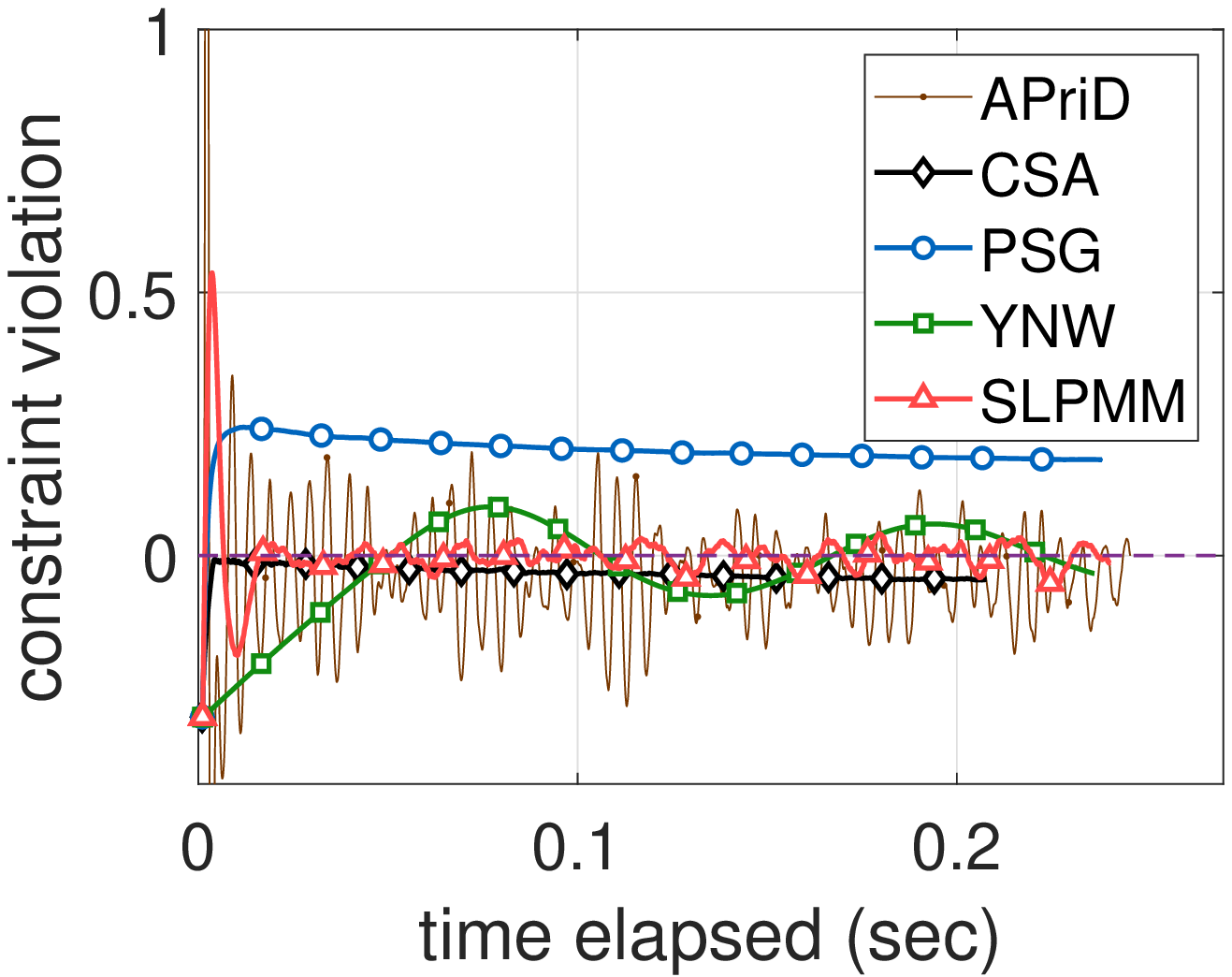}}
\end{tabular}
\caption{Comparison of algorithms on $\mathtt{CINA}$ for Neyman-Pearson classification.}
\label{figure:CINA}
\end{figure}
\begin{figure}[htp]
\centering
\setlength{\belowcaptionskip}{-6pt}
\begin{tabular}{cccc}
\subfloat[$\mathtt{objective \slash epochs}$]{
\includegraphics[width=5.5cm]{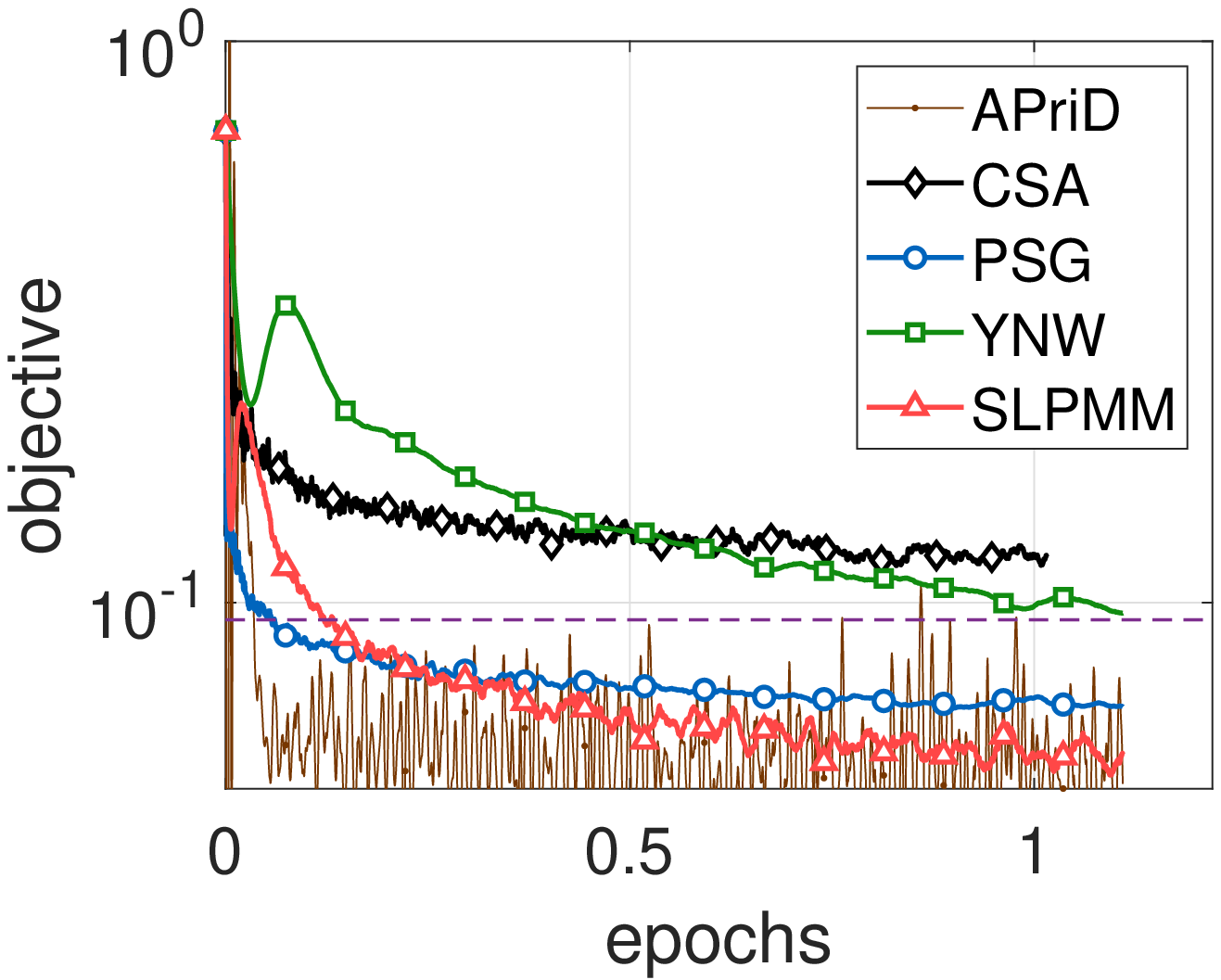}} &
\subfloat[$\mathtt{constraint \slash epochs}$]{
\includegraphics[width=5.5cm]{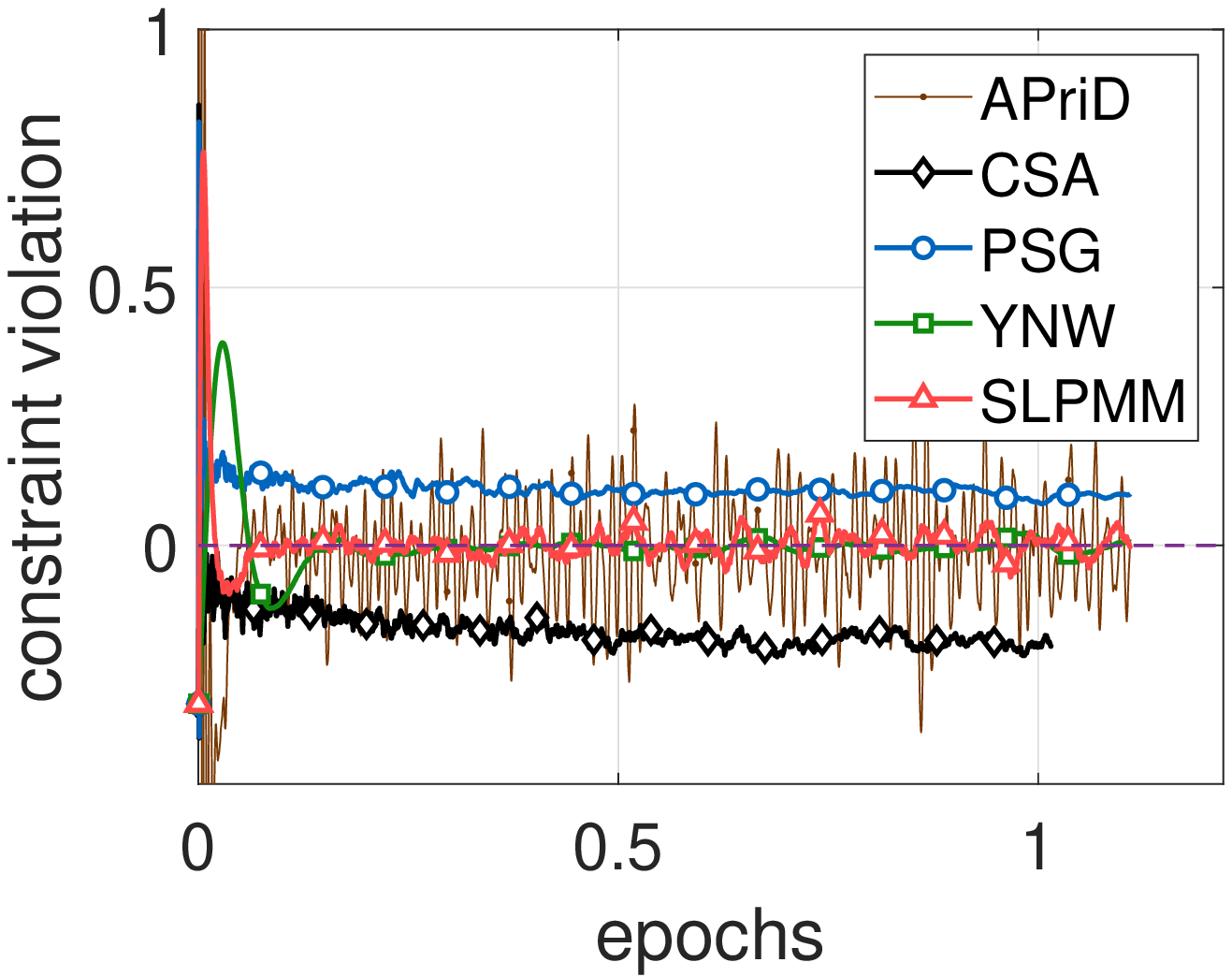}} &\\
\subfloat[$\mathtt{objective \slash cputime}$]{
\includegraphics[width=5.5cm]{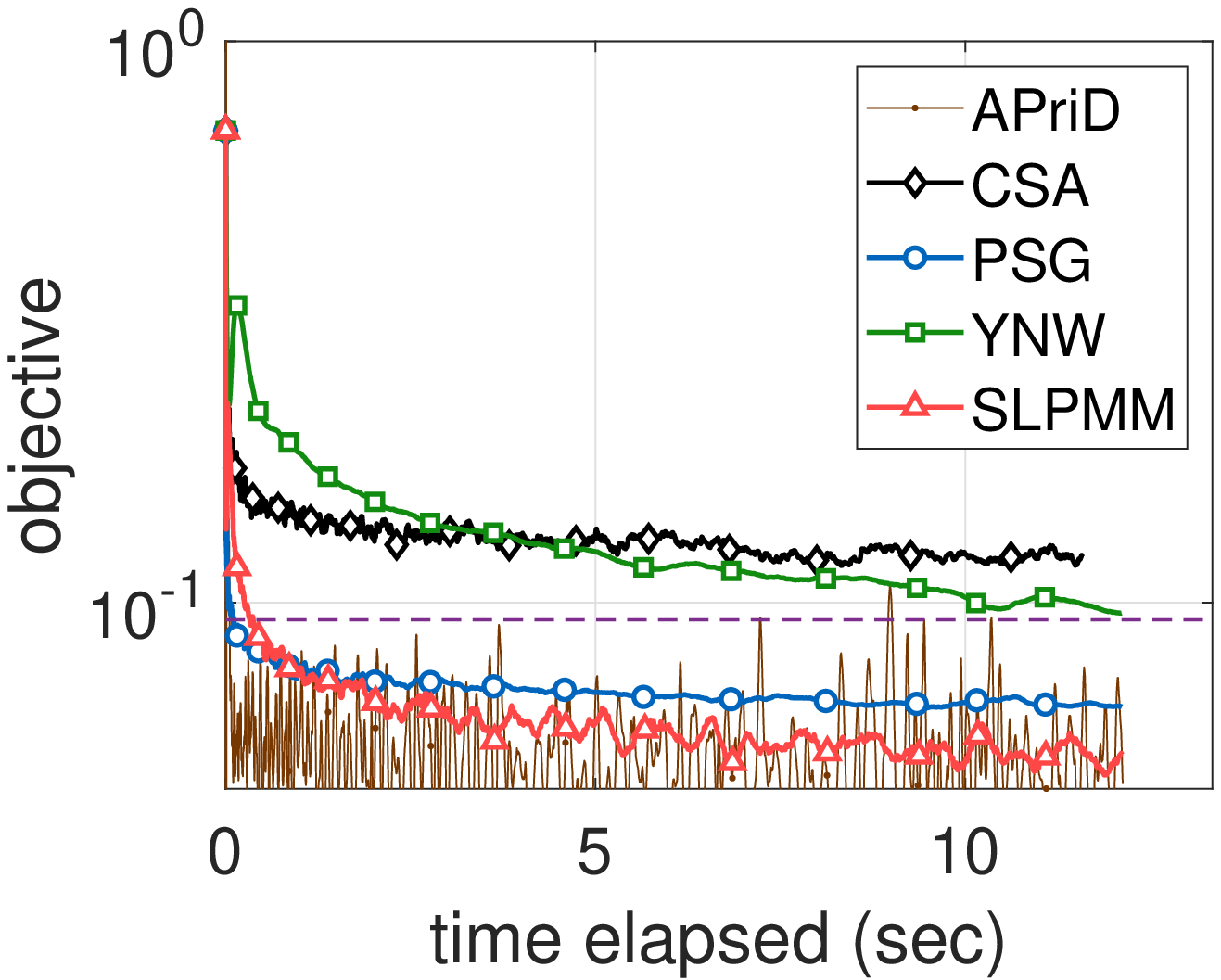}} &
\subfloat[$\mathtt{constraint \slash cputime}$]{
\includegraphics[width=5.5cm]{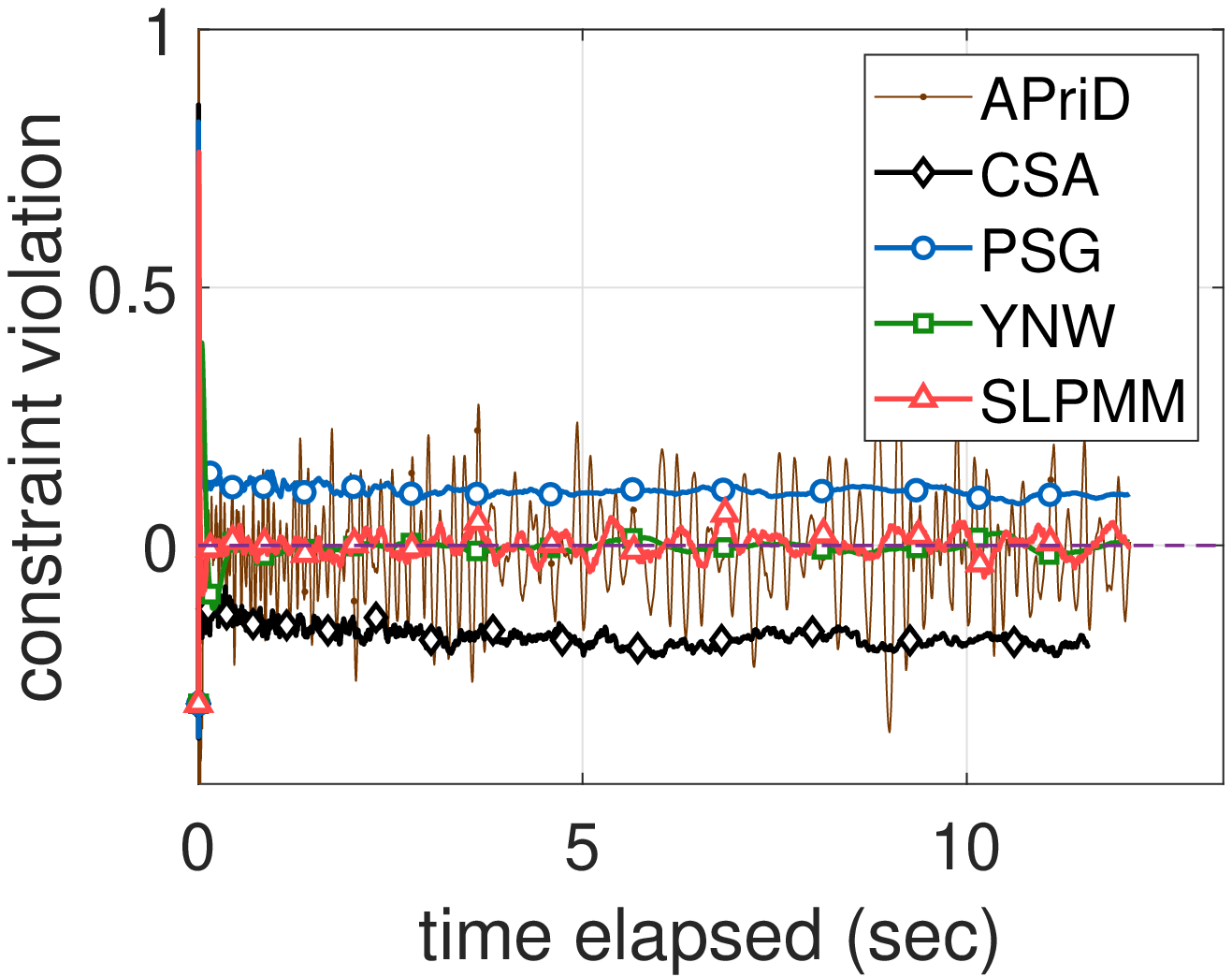}}
\end{tabular}
\caption{Comparison of algorithms on $\mathtt{MNIST}$ for Neyman-Pearson classification.}
\label{figure:MNIST}
\end{figure}
\subsection{Stochastic quadratically constrained quadratical programming}
In this subsection, we consider the following stochastic quadratically constrained quadratical programming
\[
\begin{array}{ll}
\min\limits_{x\in\cC} &f(x):=\mathbb{E}\left[\frac{1}{2}x^TA^{(0)}x+(b^{(0)})^Tx-c^{(0)}\right]\\[8pt]
\mbox{s.t.}\quad &g_i(x):=\mathbb{E}\left[\frac{1}{2}x^TA^{(i)}x+(b^{(i)})^Tx+c^{(i)}\right]\leq 0,\ i=1,2,\ldots,p,\\[5pt]
\end{array}
\]
where $A^{(i)}\in\cS_{+}^n$, $b^{(i)}\in\R^n$, $c^{(i)}\in\R$ for $i=0,1,\ldots,p$. Here, $\cS_+^n$ denotes  the set of all $n\times n$ positive semidefinite matrices. The expectations are taken with respect to the components of the parameters $\{A^{(i)},b^{(i)},c^{(i)}\}_{i=0}^p$, which are all  random variables.

The following numerical example is partially motivated by \cite{CZP2021}. The set $\cC:=\{x\in\R^n:\|x\|\leq R\}$, where $R>0$ is a constant. Let $\widehat{x}\in\R^n$ be a given point with its enty $\widehat{x}_i$ being uniformly generated from $\left(-\frac{R}{\sqrt{n}},\frac{R}{\sqrt{n}}\right)$. Let $I_n$ be the identity matrix. For each $i=0,1,\ldots,p$, the random matrix $A^{(i)}=I_n+\Delta_i$, where $\Delta_i$ is a symmetric matrix and its entry is uniformly distributed over $[-0.1,0.1]$. The random vector $b^{(i)}$ is uniformly distributed from $[-1,1]$. The random variable $c^{(i)}$ is constructed with a particular purpose. Let $h^{(i)}$ be a random variable uniformly distributed over $[0,2i]$, then define $c^{(i)}=-(\frac{1}{2}\widehat{x}^TA^{(i)}\widehat{x}+(b^{(i)})^T\widehat{x}+h^{(i)})$. In this setting, we can easily verify that $g_i(\widehat{x})=-i<0$ for $i=1,\ldots,p$ and hence the Slater's condition is satisfied. We can also get that the optimal solution is $0$ and the optimal value is $\frac{1}{2}\|\widehat{x}\|^2$.

In this experiment, we compare the performance of SLPMM with PSG, YNW and APriD. At each iteration of the algorithms, we generate the samples of  $\{A^{(i)},b^{(i)},c^{(i)}\}_{i=0}^p$ based on the above distributions for function and gradient evaluation. We set $n=100$, $p=5$, $R=2$. The maximum number of iterations is set to $K=1000$. The initial point is set to $x^0=(\sqrt{R/n},\sqrt{R/n},\ldots,\sqrt{R/n})^T$.

The results in terms of time are shown in Figure \ref{figure:QCQP}. From picture (b) (plots the value of $\max_i\{g_i(x^k)\}$), we can see the iterations of all algorithms satisfy the constraints. From picture (a), we observe that SLPMM is comparable with PSG, and obviously outperforms over APriD and YNW.
\begin{figure}[htp]
\centering
\setlength{\belowcaptionskip}{-6pt}
\begin{tabular}{cccc}
\subfloat[$\mathtt{objective \slash cputime}$]{
\includegraphics[width=5.5cm]{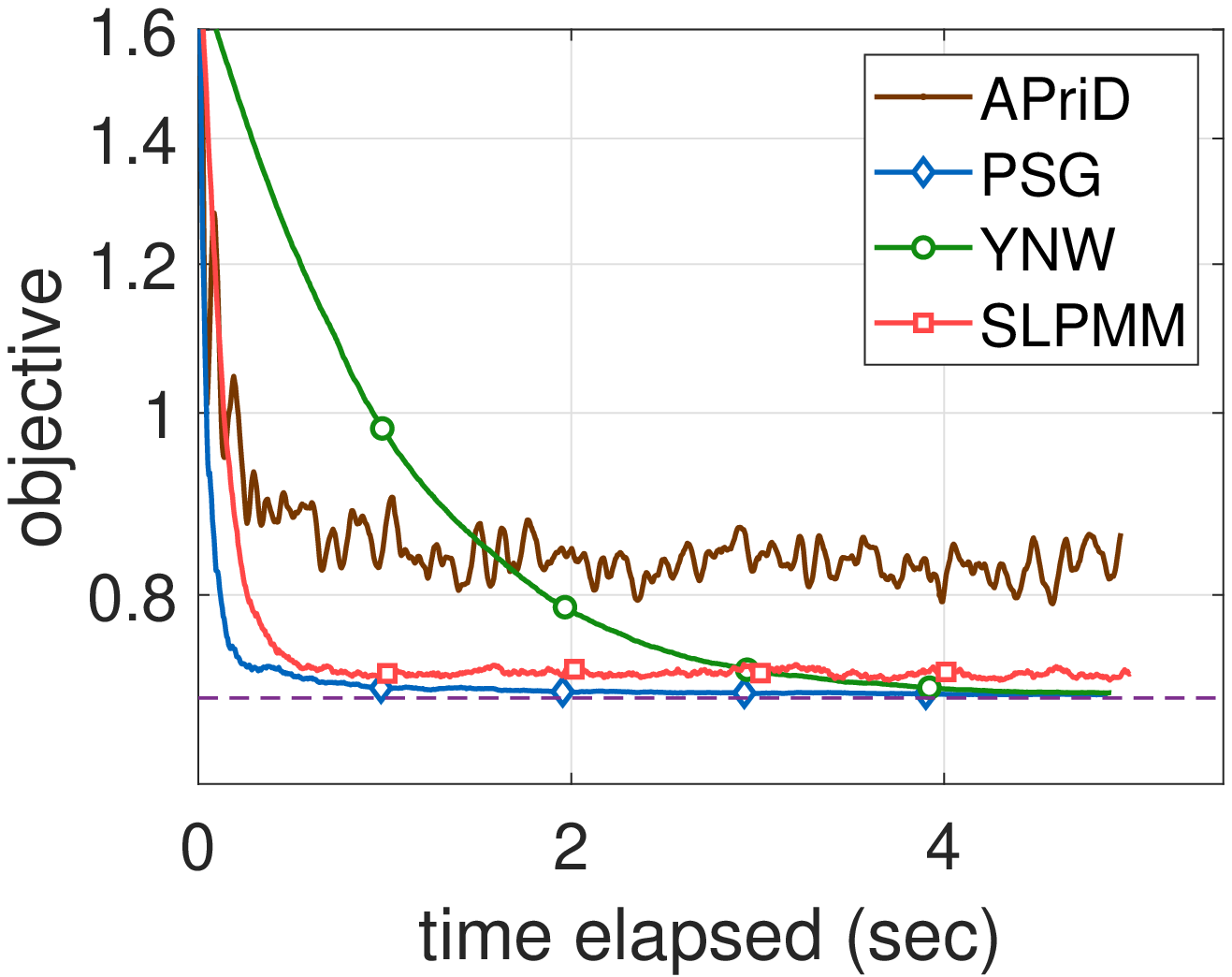}} &
\subfloat[$\mathtt{constraint \slash cputime}$]{
\includegraphics[width=5.5cm]{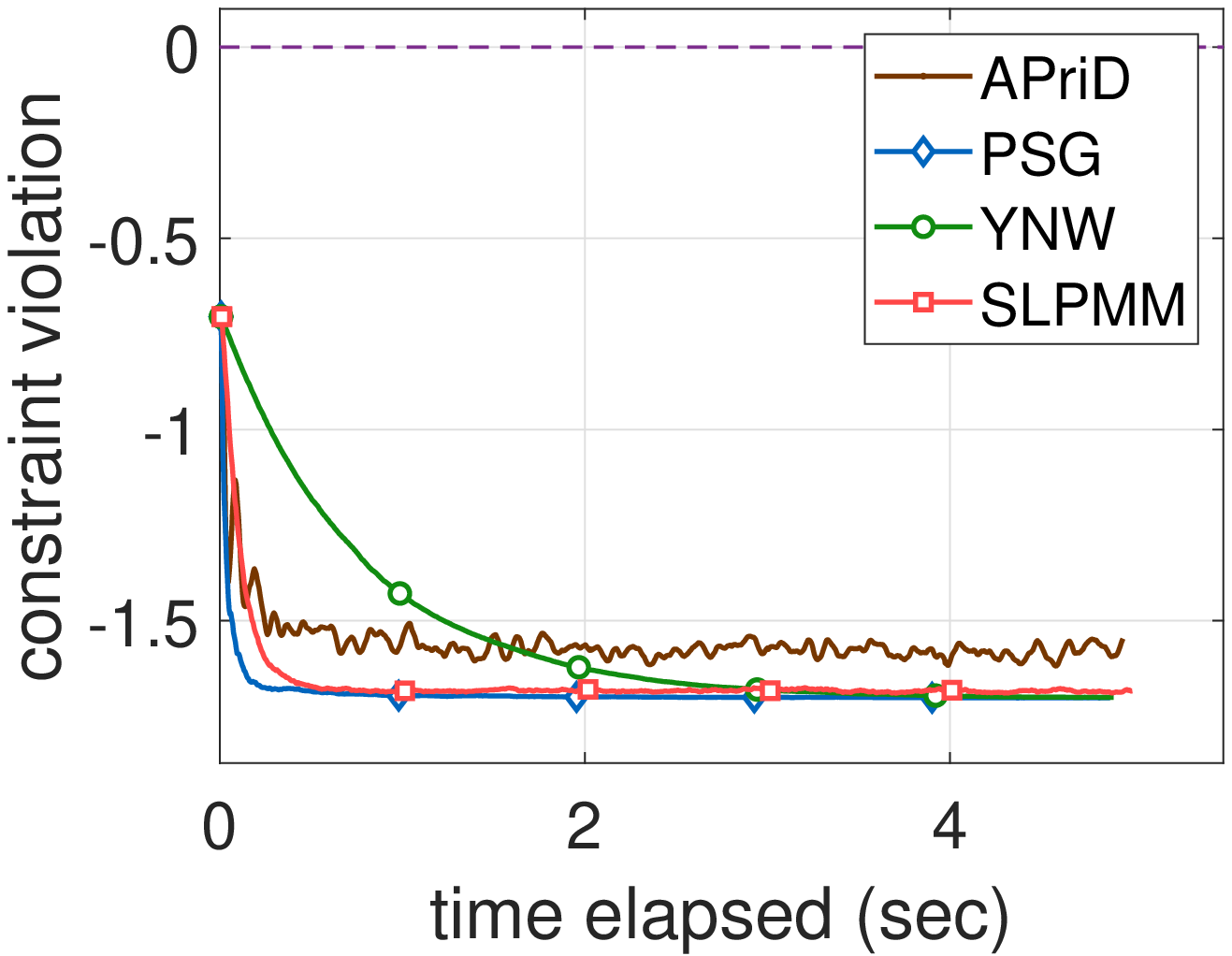}}
\end{tabular}
\caption{Comparison of algorithms on stochastic quadratically constrained quadratical programming.}
\label{figure:QCQP}
\end{figure}
\subsection{Second-order stochastic dominance  constrained portfolio optimization}
In this subsection, we consider the following second-order stochastic dominance (SSD) constrained portfolio optimization problem
\[
\begin{array}{ll}
\min &\mathbb{E}[-\xi^Tx]\\[5pt]
\mbox{s.t.}\quad &\mathbb{E}[[\eta-\xi^Tx]_+]\leq\mathbb{E}[[\eta-Y]_+],\quad\forall \eta\in\R,\\[5pt]
&x\in\cC:=\{x\in\R^n:\sum_{i=1}^nx_i=1,\ \bar{x}\geq x\geq 0\},
\end{array}
\]
where $\bar{x}$ is the upper bound and $Y$ stands for the random return of a benchmark portfolio dominated by the target portfolio in the SSD sense. Since it was first introduced by \cite{DR2003}, SSD has been widely used to control risk in financial portfolio  \citep{KD2018,Noyan2018}.
\cite{KKU2016} showed that,
if $Y$ is  discretely distributed with $\{y_1,y_2,\ldots,y_p\}$, the SSD constrained portfolio optimization is reduced to
\begin{equation}\label{eq:port-ssd}
\begin{array}{ll}
\min &f(x):=\mathbb{E}[-\xi^Tx]\\[8pt]
\mbox{s.t.}\quad &g_i(x):=\mathbb{E}[[y_i-\xi^Tx]_+]-\mathbb{E}[[y_i-Y]_+]\leq 0,\quad i=1,\ldots,p,\\[8pt]
&x\in \cC:=\{x\in\R^n:\sum_{i=1}^nx_i=1,\ \bar{x}\geq x\geq 0\},
\end{array}
\end{equation}
which is  an instance of Problem (\ref{eq:1}).

\cite{DMW2016} proposed several methods for solving SSD constrained optimization problems based on augmented Lagrangian framework and analyze their convergence. In particular, the proposed approximate augmented Lagrangian method with exact minimization (PALEM) has some similarities to SLPMM. At each iteration in PALEM, a minimization problem with respect to the augmented Lagrangian function of a reduced problem is solved to obtain $x^k$, and the multiplier $\mu^k$ is updated. They proved that the sequences $\{x^k\}$ and $\{\mu^k\}$ converge to the optimal solution of primal and dual problem, respectively. In contrast, although SLPMM is  also constructed based on augmented Lagrangian framework as PALEM, they are quite different. The subproblem at each iteration in SLPMM is a minimization problem of a linearized augmented Lagrangian function together with a proximal term, which is easier to solve. The sampling strategy is different. In PALEM, the sample set is updated at each iteration based on the calculation of the expectation of constraint function. SLPMM only simply requires one sample at each iteration.  Moreover, since in our setting the expectation is assumed to be impossible to be calculated, we can not obtain the convergence of the sequence to optimal solution.

In this experiment, we compare the performance of SLPMM with APriD, PSG, YNW and PALEM to solve Problem (\ref{eq:port-ssd}) on  the following four datasets
\[
\{``\texttt{Dax\_26\_3046}",``\texttt{DowJones\_29\_3020}",``\texttt{SP100\_90\_3020}",``\texttt{DowJones\_76\_30000}"\}
\] from \citep{KKU2016}. Take ``\texttt{DowJones\_29\_3020}" for example, ``\texttt{DowJones}" stands for Dow Jones Index, 29 is the number of stocks  and 3020 is the number of scenarios, i.e.,  $n=29, p=3020$. The initial point  is set to 0. For PALEM, we use the MATLAB function \texttt{fmincon} to solve the subproblem. For SLPMM, we utilize the Nesterov's accelerated projected gradient method (APG) to solve the subproblem (\ref{eq:general-subp}), the stopping criterion of APG is set to $\|y^t-T_{L_t}(y^t)\|\leq 10^{-6}$, and the projection $\Pi_{\cC}$ is computed by the method proposed in \citep{WL2015}. In particular,    since the number of the constraints of Problem (\ref{eq:port-ssd}) is large, we apply a sampling technique to reduce the computational cost. In specific, at each iteration, instead of using the whole constraint index set $\{1,\ldots,p\}$ in the augmented Lagrangian function (\ref{augL}), we first randomly sample a subset $I_k\subset\{1,\ldots,p\}$ and then replace $\sum_{i=1}^p$ with $\sum_{i\in I_k}$ in (\ref{augL}). This sampling strategy, which is also used in \citep{Xiao2019}, is proven to be very efficient in practice. Let us also remark that, by taking an extra expectation with respect to $I_k$, the expected convergence rates of SLPMM coupled with this sampling strategy can be established in a similar way as in Section \ref{sec:rates}. This is also pointed out in \citep[Section 5]{Xiao2019}.

The numerical results are presented in Figure \ref{figure:ssd}. Since the maximum of $p$ constraint values are always zero (which indicates that the constraints are satisfied), we omit the presentation of constraint violation.  We only report the change of the objective value with respect to cputime.  The horizontal dashed line in each picture represents a reference optimal objective value which is obtained from \citep{KKU2016}. In general, we can observe that SLPMM has an obvious advantage compared with the other four algorithms. In view of dataset ``\texttt{DowJones\_76\_30000}" which refers to a large scale optimization problem with 30,000 constraints, SLPMM converges to the optimal objective value less than 4 seconds. We can also observe that SLPMM is very robust and stable for all four datasets.
\begin{figure}[htp]
\centering
\setlength{\belowcaptionskip}{-6pt}
\begin{tabular}{cccc}
\subfloat[\texttt{Dax\_26\_3046}]{
\includegraphics[width=5.5cm]{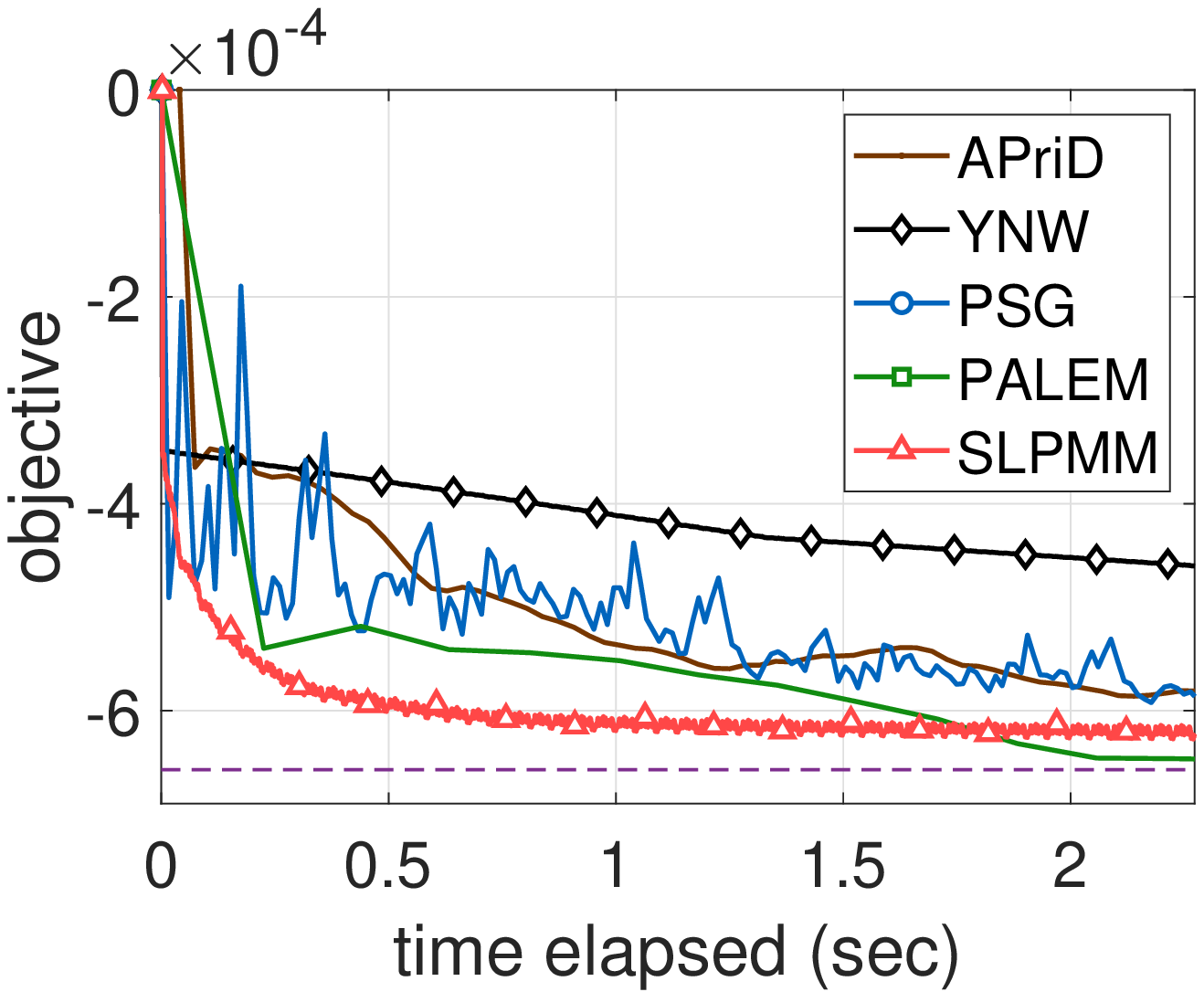}} &
\subfloat[\texttt{DowJones\_29\_3020}]{
\includegraphics[width=5.5cm]{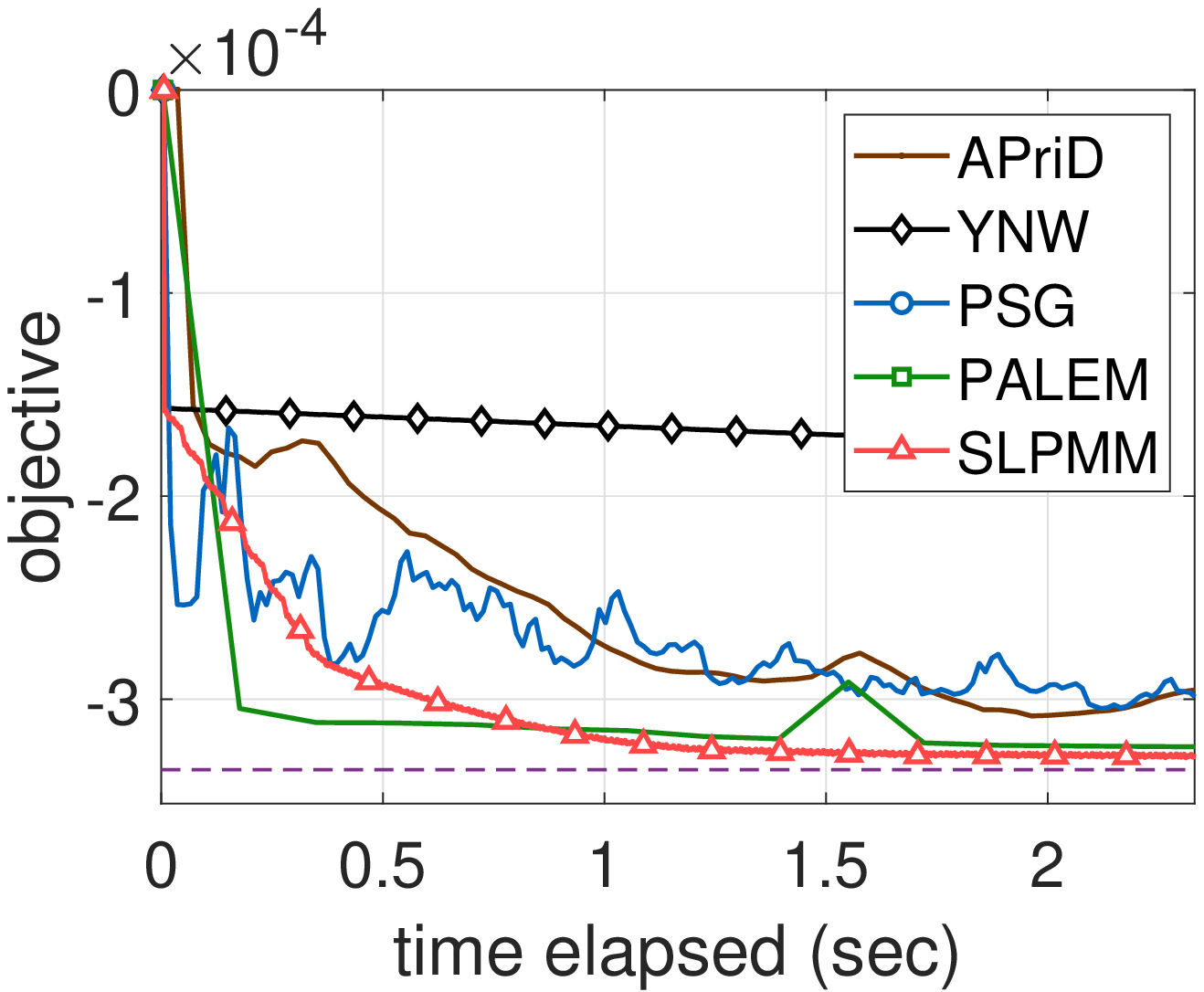}} &\\
\subfloat[\texttt{SP100\_90\_3020}]{
\includegraphics[width=5.5cm]{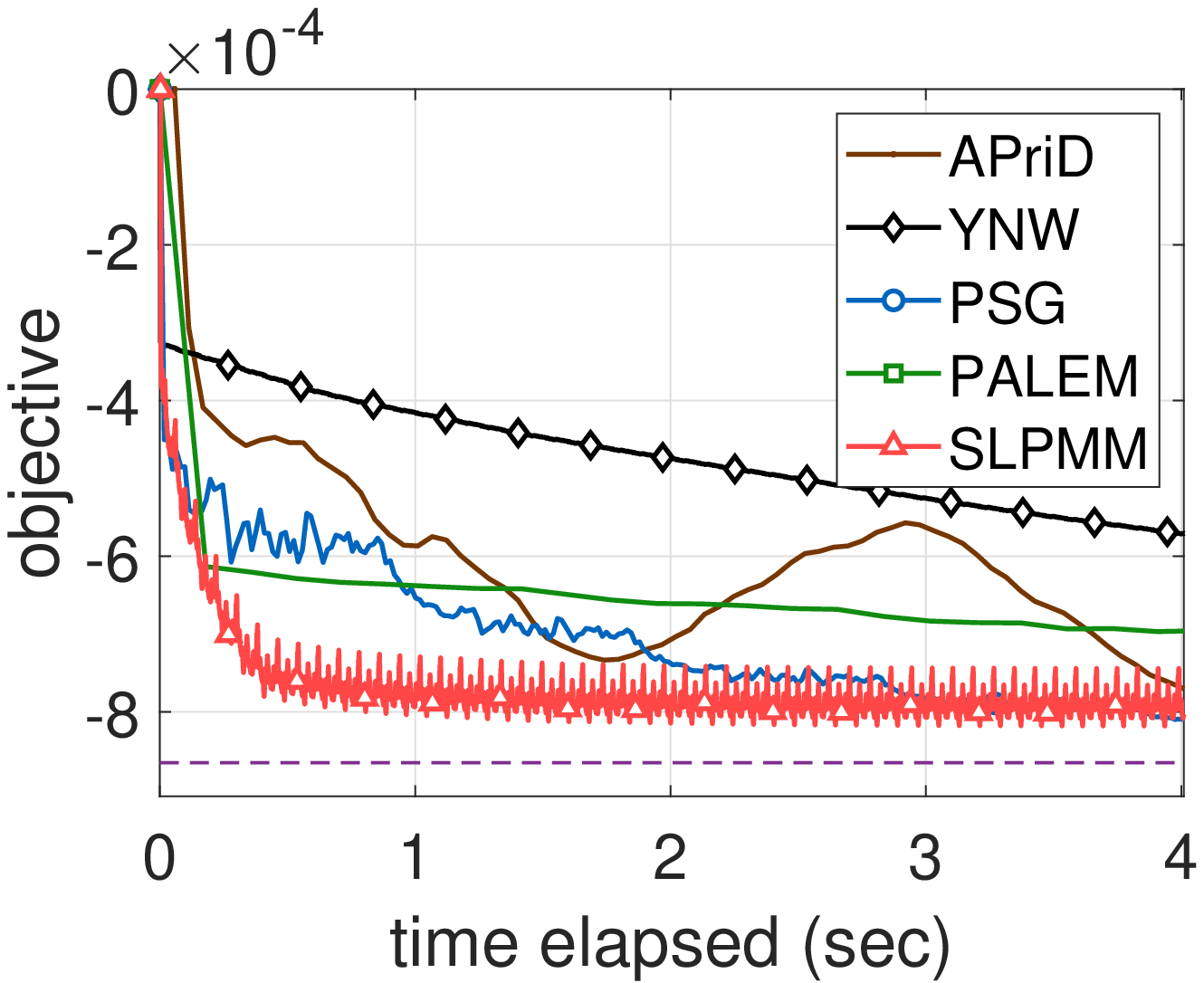}} &
\subfloat[\texttt{DowJones\_76\_30000}]{
\includegraphics[width=5.5cm]{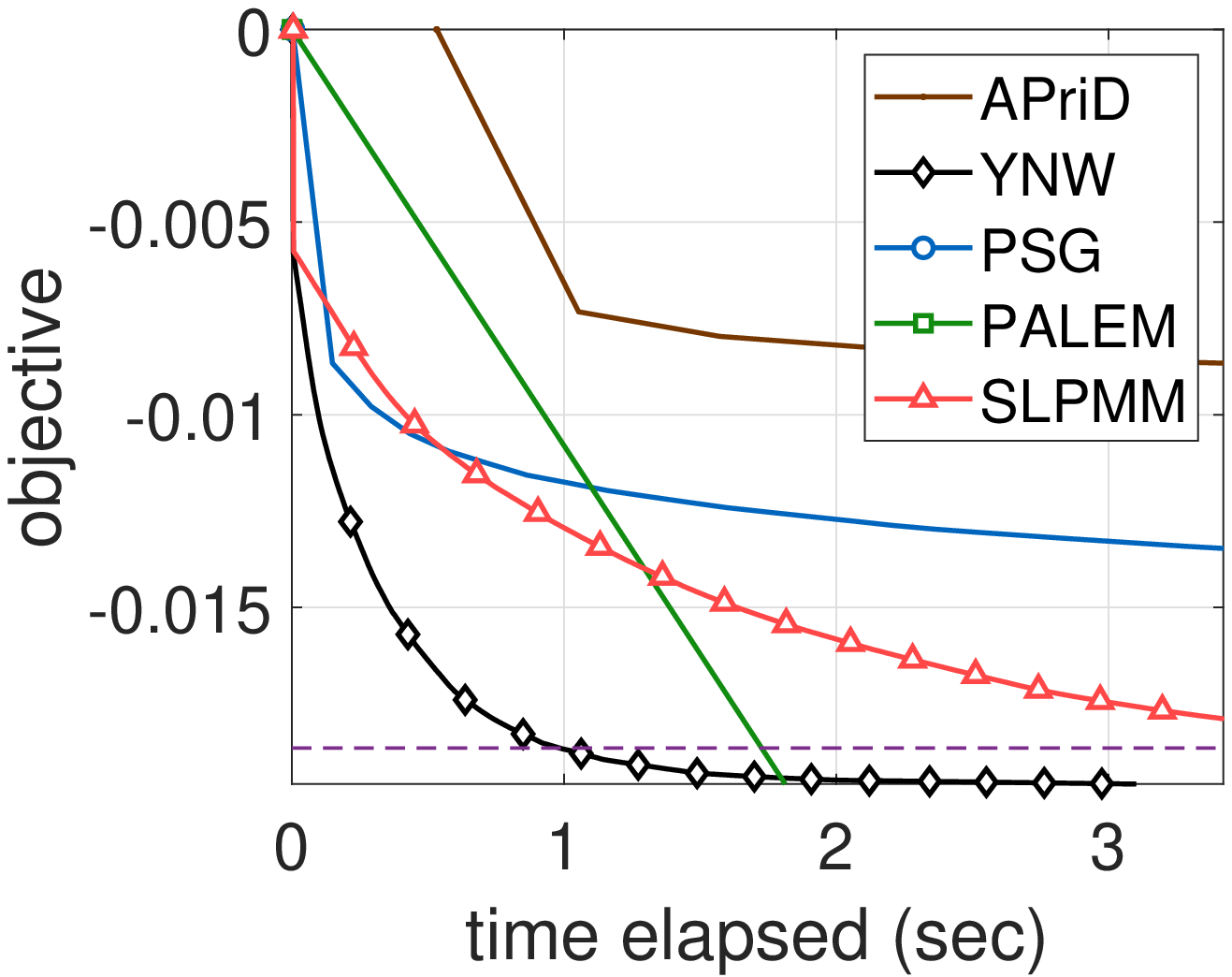}}
\end{tabular}
\caption{Comparison of algorithms for SSD constrained portfolio optimization.}
\label{figure:ssd}
\end{figure}

\section{Conclusion}\label{sec:conclusion}
We present a hybrid method of stochastic approximation technique and proximal augmented Lagrangian method. It is  shown that the expected convergence rates and the large-deviation properties  are comparable with the existing related stochastic methods. On the other hand,  the proposed method is parametric-independent.  Numerical experiments also demonstrate the superiority  in comparison with  the  stochastic first-order methods. Thus, both theoretical and numerical results suggest that the proposed algorithm is efficient for solving convex stochastic programming with expectation constraints.

However, there are still several valuable questions left to be answered. It is well-known that the deterministic augmented Lagrangian can achieve superlinear convergence. Therefore, the first question is whether the convergence rates can be improved to match the numerical performance and the rates in the deterministic setting. Secondly, it is worthwhile to consider the inexact method, that is, the subproblem is solved inexactly.
Another interesting  topic is how to use the techniques in this paper to deal with  nonconvex stochastic optimization. The proposed algorithm in the current form is not applicable to solve nonconvex problems, such as chance constrained programs \citep{BSZ2021} and MIMO transmit signal design problem \citep{LLK2019}.

Finally, let us mention that the stochastic algorithms for stochastic optimization can be easily extended to solve online problems, and vice versa, see \citep{YMNeely2017} for instance. Hence, the proposed method can be slightly revised to solve the corresponding online problems.

%


\ACKNOWLEDGMENT{%
 The authors would like to thank the anonymous reviewers and the associate
editor for the valuable comments and suggestions that
helped us to greatly improve the quality of the paper.
}

%
%
%


\bibliographystyle{informs2014} 
\bibliography{ref.bib} 



\end{document}